\RequirePackage{etoolbox}
\csdef{input@path}{{style/}{graphics/}}
\documentclass[aos,MSNbibl,nameyear,dvips]{arximspdf}
\usepackage{mathbh}
\usepackage{accents}
\usepackage{graphicx}

\doi{10.1214/14-AOS1305}
\volume{43}
\issue{2}
\pubyear{2015}
\firstpage{903}
\lastpage{934}
\docsubty{FLA}

\makeatletter
\renewcommand{\mathring}[1]{\accentset{\circ}{#1}}
\newcommand{\rrvert}{\vert}
\newcommand{\rrVert}{\Vert}
\newcommand{\llvert}{\vert}
\newcommand{\llVert}{\Vert}
\newcommand{\iint}{\int\!\!\!\int}
\newcommand{\overset}{\stackrel}
\def\pr{\mathbb P}
\def\ind{\mathbh{1}}
\newtheorem{theorem}{Theorem}
\newtheorem{lemma}[theorem]{Lemma}
\newtheorem{proposition}[theorem]{Proposition}
\newproclaim{remark}{Remark}
\makeatother

\begin{document}
\begin{frontmatter}

\title{Bias correction in multivariate extremes\thanksref{T1}}
\runtitle{Bias correction in multivariate extremes}

\begin{aug}
\author[A]{\fnms{Anne-Laure}~\snm{Foug\`eres}\corref{}\thanksref{M1}\ead[label=e1]{fougeres@math.univ-lyon1.fr}},
\author[B]{\fnms{Laurens}~\snm{de Haan}\thanksref{M2}\ead[label=e2]{ldehaan@ese.eur.nl}}\\
\and
\author[A]{\fnms{C\'ecile}~\snm{Mercadier}\thanksref{M1}\ead[label=e3]{mercadier@math.univ-lyon1.fr}}
\runauthor{A.-L. Foug\`eres, L. de Haan and C. Mercadier}
\affiliation{Universit\'e Lyon 1\thanksmark{M1} and Erasmus
University\thanksmark{M2}}
\address[A]{A.-L. Foug\`eres\\
C. Mercadier\\
Universit\'e de Lyon, CNRS, Universit\'e Lyon 1\\
Institut Camille Jordan\\
43 blvd du 11 novembre 1918\\
F-69622 Villeurbanne-Cedex\\
France\\
\printead{e1}\\
\phantom{E-mail: }\printead*{e3}}
\address[B]{L. de Haan\\
Department of Economics\\
Erasmus University\\
P.O. Box 1738\\
3000 DR Rotterdam\\
The Netherlands\\
\printead{e2}}
\end{aug}
\thankstext{T1}{Supported in part by the Agence Nationale de la Recherche
through the AST\&RISK project (ANR-08-BLAN-0314-01), and by
FCT/PTDC/MAT/112770/2009 (Portugal).}

%
\received{\smonth{2} \syear{2014}}
%
\revised{\smonth{12} \syear{2014}}

%
\begin{abstract}
The estimation of the extremal dependence structure is spoiled by the
impact of the bias, which increases
with the number of observations used for the estimation. Already known
in the univariate setting,
the bias correction procedure is studied in this paper under the
multivariate framework.
New families of estimators of the stable tail dependence function are
obtained. They are
asymptotically unbiased versions of the empirical estimator introduced
by Huang [Statistics of bivariate extremes (1992) Erasmus Univ.].
Since the new estimators have a regular behavior with respect to the
number of observations,
it is possible to deduce aggregated versions so that the choice of the
threshold is substantially
simplified. An extensive simulation study is provided as well as an
application on real data.
\end{abstract}

%
\end{frontmatter}

\section{Introduction}
Estimating extreme risks in a multivariate framework is highly
connected with the estimation of the extremal dependence structure.
This structure can be described \textit{via} the stable tail
dependence function (s.t.d.f.) $L$, first introduced by \citet{Huang1992}.
For any arbitrary dimension $d$, consider a multivariate vector
$(X^{(1)},\ldots,X^{(d)})$ with continuous marginal cumulative
distribution functions (c.d.f.) $F_1, \ldots, F_d$. The s.t.d.f. is defined
for each positive reals $x_1,\ldots,x_d$ as
\begin{eqnarray*}
&& \lim_{t \to\infty} t \pr\bigl\{ 1-F_1\bigl(X^{(1)}
\bigr)\leq t^{-1}x_1\mbox{ or } \ldots\mbox{ or }
1-F_d\bigl(X^{(d)}\bigr) \leq t^{-1}x_d
\bigr\}
\nonumber\\[-8pt]\\[-8pt]\nonumber
&&\qquad = L(x_1,\ldots,x_d).
\end{eqnarray*}
Assuming that such a limit exists and is nondegenerate is equivalent to
the classical assumption of
existence of a multivariate domain of attraction for the componentwise
maxima; see, for example,
\citet{dehaanferreira2006}, Chapter~7.
The previous limit can be rewritten as
%
\begin{equation}
\label{eqL} \lim_{t \to\infty} t \bigl[ 1- F\bigl
\{F_1^{-1}\bigl(1-t^{-1}x_1\bigr),
\ldots, F_d^{-1}\bigl(1-t^{-1}x_d
\bigr) \bigr\}\bigr] = L(x_1,\ldots,x_d),
\end{equation}
where $F$ denotes the multivariate c.d.f. of the vector $(X^{(1)},\ldots,X^{(d)})$,
and for $j=1,\dots,d$, $F_j^{-1}(t) = \inf\{z \in{\mathbb R}\dvtx  F_j(z)
\geq t \}$.
Consider\vspace*{1pt} a sample of size $n$ drawn from~$F$ and an intermediate
sequence, that is to say a sequence $k=k(n)$ tending to infinity as $n
\to\infty$, with $k/n \to0$. Denote by ${\mathbf{x}}=(x_1,\ldots,x_d)$
a vector of the positive quadrant $\mathbb{R}^d_+=\{ (x_1,\ldots,x_d)\dvtx  x_j \geq0, 
j=1,\dots, d\}$ and by $X^{(j)}_{k,n}$ the $k$th order statistics
among $n$ realizations of the margins $X^{(j)}$.
The empirical estimator of $L({\mathbf{x}})$ is obtained from~(\ref
{eqL}), replacing $F$ by its empirical version, $t$ by $n/k$, and
$F_j^{-1}(1-t^{-1}x_j)$ for $ j=1,\ldots,d$ by its empirical counterpart
$X^{(j)}_{n-[nt^{-1}x_j],n}$, so that
%
\begin{equation}
\label{eq1storderhat} \hat L_k({\mathbf x})= \frac{1}k \sum
_{i=1}^n \ind_{ \{X^{(1)}_i \geq
X^{(1)}_{n-[kx_1]+1,n}~\mathrm{or}~\ldots~\mathrm{or}~X^{(d)}_i \geq X^{(d)}_{n-[kx_d]+1,n}  \} }.
\end{equation}
See \citet{Huang1992} for pioneering works on this estimator. Under
suitable conditions, it can be shown (see Section~\ref{secnotation})
that the estimator $\hat L_k({\mathbf x})$ has the following asymptotic expansion:
%
\begin{equation}
\label{eqasympt-expL} \hat L_k({\mathbf x}) - L({\mathbf x}) \approx
\frac{Z_L({\mathbf x})}{\sqrt{k}} + \alpha(n/k)M({\mathbf x}),
\end{equation}
where $Z_L$ is a continuous centered Gaussian process, $\alpha$ is a
function that tends to~0 at infinity and $M$ is a continuous function.
In particular\break $ \sqrt{k} \{ \hat L_k({\mathbf x}) - L({\mathbf x}) \}$ can be
approximated in distribution by $Z_L({\mathbf x})$, provided that $ \sqrt
{k} \alpha(n/k)$ tends to 0 as $n$ tends to infinity.
This condition imposes a slow rate of convergence of the estimator
$\hat L_k({\mathbf x})$, so one would be interested in relaxing this
hypothesis. As a counterpart, as soon as $ \sqrt{k} \alpha(n/k)$
tends to a nonnull constant~$\lambda$, an asymptotic bias appears and
is explicitely given by $\lambda M({\mathbf x})$. The aim of this paper is
to provide a procedure that reduces the asymptotic bias. The latter
will be estimated and then subtracted from the empirical estimator.
This kind of approach has been considered in the univariate setting for
the bias correction of the extreme value index with unknown sign by
\citet{caidehaanzhou2013}. Refer also to \citeauthor{peng1998} (\citeyear{peng1998,peng2010}) %
\citet{fragadehaanlin2003}, \citet{gomesdehaanrodrigues2008} and
\citet{caeirogomesrodrigues2009} for previous
contributions on this problem.
Note finally that the case of dependent sequences has been recently
studied by \citet{dehaanmercadierzhou2014}.

The nonparametric estimation of the extremal dependence structure has
been widely studied in the bivariate case; see, for instance, \citet
{Huang1992}, \citet{einmahldehaansinha1997}, \citet
{caperaafougeres2000}, \citet{abdousghoudi2005}, \citet
{guillotteperronsegers2011} and \citet{bucherdettevolgushev2011}.
Bias correction problems in the bivariate context received less
attention than in the univariate setting. To the best of our knowledge,
it seems to be reduced to \citet{beirlantdierckxguillou2011} and
\citet{goegebeurguillou2013}, who consider the estimation of
bivariate joint tails, which differs slightly from our task.

As for the multivariate framework, \citet{dehaanresnick1993}
introduces the empirical estimator. General approaches under parametric
assumptions on the function $L$ have been developed, for example, by
\citet{colestawn1991}, \citet{joesmithweissman1992}, \citet
{einmahlkrajinasegers2008} and \citet{einmahlkrajinasegers2012}.
Apparently, no procedure correcting the bias can be found in the
literature for dimension greater than two. The objective of this
article is to fill this gap. Note that our method does not only consist
of applying the univariate bias procedure at several points. Indeed,
the bias is no longer a parametric function, so that the new feature is
mainly the fact that we are able to estimate and then subtract a
function with an unknown form.
Two families of asymptotically unbiased estimators of the s.t.d.f. are
proposed, and their theoretical behaviors are studied. A practical
advantage of these new estimators is that they can be aggregated, thus
reducing the variability.

The paper is organized as follows: Section~\ref{secnotation} contains
hypotheses and first results. The bias reduction procedure is described
in Section~\ref{secprocedure}, and the main theoretical results are
presented therein. Several theoretical models are exhibited in
Section~\ref{secexamples} that satisfy the required assumptions.
Section~\ref{secsimu} illustrates the performance of the new
estimators on both simulated and real data. The estimation of side
components is postponed up to Section~\ref{secrho}. The proofs are
relegated to Section~\ref{secproofs}.

\section{Notation, assumptions and first results} \label{secnotation}
Let ${\mathbf X}_1=(X^{(1)}_1,\ldots,X^{(d)}_1),\break \ldots, {\mathbf
X}_n=(X^{(1)}_n, \ldots,X^{(d)}_n)$ be independent and identically
distributed multivariate random vectors with c.d.f. $F$ and continuous
marginal c.d.f.'s $F_j$ for $j=1,\ldots,d$. Suppose $F$ is in the domain
of attraction of an extreme value distribution with c.d.f.~$G$. We recall
that it supposes the existence for $j=1,\ldots,d$ of sequences
$a_n^{(j)}>0$, $b_n^{(j)}$ of real numbers and a c.d.f. $G$ with
nondegenerate marginals such that
\begin{eqnarray*}
&& \lim_{n\to\infty} \mathbb{P}\bigl(\max\bigl\{X^{(1)}_1,
\ldots,X^{(1)}_n\bigr\} \leq a^{(1)}_n
x_1 + b^{(1)}_n, \ldots,
\\
&&\hspace*{54pt} \max\bigl
\{X^{(d)}_1,\ldots,X^{(d)}_n\bigr\}\leq a^{(d)}_n x_d + b^{(d)}_n \bigr)=G({\mathbf x})
\end{eqnarray*}
for all points ${\mathbf x}$ where $G$ is continuous. Denote by $G_j$ the
$j$th marginal c.d.f. of $G$. It is possible to show that the domain of
attraction condition can be expressed as condition~(\ref{eqL}) along
with the convergence of the marginal distributions to the $G_j$'s, and that
%
\begin{equation}
\label{eqLbis} L({\mathbf x})=-\log G \bigl(\{-\log G_1\}
^{-1}(x_1),\ldots,\{-\log G_d
\}^{-1}(x_d) \bigr).
\end{equation}
Let $\mu$ be the measure defined by
%
\begin{equation}
\label{eqmu} \mu\bigl\{A({\mathbf x})\bigr\}:=L({\mathbf x}),
\end{equation}
where $A({\mathbf x}):=\{{\mathbf u}\in\mathbb{R}_+^d$: there exists
$j$ such that $u_j>x_j \}$ for any vector ${\mathbf{x}}\in \mathbb{R}_+^d$.

Several conditions are now described. The first two have been
introduced by \citet{dehaanresnick1993}:
\begin{itemize}[--]
\item[--] The \textit{first-order condition} consists of assuming that the
limit given in~(\ref{eqL}) exists, and that the convergence is
uniform on any $[0,T]^d$, for $T>0$.
\item[--] The \textit{second-order condition} consists of assuming the
existence of a positive 
function $\alpha$, such that $\alpha(t) \to0$ as $t\to\infty$, and
a nonnull function $M$ such that for all ${\mathbf x}$ with positive coordinates,
%
\begin{eqnarray}\label{eq2ndorder}
&& \lim_{t \to\infty} \frac{1}{\alpha(t)} \bigl\{ t \bigl[ 1-
F\bigl\{ F_1^{-1}\bigl(1-t^{-1}x_1
\bigr),\ldots, F_d^{-1}\bigl(1-t^{-1}x_d
\bigr) \bigr\}\bigr] - L({\mathbf x}) \bigr\}
\nonumber\\[-8pt]\\[-8pt]\nonumber
&&\qquad = M({\mathbf x}),
\end{eqnarray}
uniformly on any $[0,T]^d$, for $T>0$.
\item[--] The \textit{third-order condition} consists of assuming the
existence of a positive 
function $\beta$, such that $\beta(t) \to0$ as $t\to\infty$, and a
nonnull function $N$ such that for all ${\mathbf x}$ with positive coordinates,
%
\begin{eqnarray}\label{eq3rdorder}
&& \lim_{t \to\infty} \frac{1}{\beta(t)} \biggl\{\frac{t   [ 1- F\{
F_1^{-1}(1-t^{-1}x_1),\ldots,   F_d^{-1}(1-t^{-1}x_d) \}] - L({\mathbf
x})}{\alpha(t)} - M({\mathbf x}) \biggr\}
\nonumber\\[-8pt]\\[-8pt]\nonumber
&&\qquad  = N({\mathbf x}),\hspace*{-10pt}
\end{eqnarray}
uniformly on any $[0,T]^d$, for $T>0$. This implicitly requires that
$N$ is not a multiple of the function $M$; see Remark \ref{rmMN}.
\end{itemize}

%
\begin{remark}\label{rkhomoL}
The function $L$ defined by~(\ref
{eqL}) and that appears in~(\ref{eq2ndorder}) and~(\ref
{eq3rdorder}) is homogeneous of order 1. We refer, for instance, to
\citet{dehaanferreira2006}, pages 213~and~236. Most of the
estimators constructed in this paper use the homogeneity property. Note\vspace*{1pt}
that pointwise convergence in~(\ref{eqL}) entails uniform convergence
on the square $[0,T]^d$. See, for instance, \citet{dehaanferreira2006}, page~237.
\end{remark}

%
\begin{remark}\label{rmMN} If $N=c\cdot M$ for some constant $c$,
the relation can be reformulated as
\begin{eqnarray*}
&& \lim_{t \to\infty} \frac{1}{\beta(t)} \biggl\{ \frac{t   [ 1- F\{
F_1^{-1}(1-t^{-1}x_1),\ldots,   F_d^{-1}(1-t^{-1}x_d) \}] - L({\mathbf
x}) }{\alpha(t)(1+c\beta(t))} - M({
\mathbf x}) \biggr\}
\\
&&\qquad = 0,
\end{eqnarray*}
which we want to exclude. We refer to \citeauthor{dehaanferreira2006} [(\citeyear{dehaanferreira2006}), page~385],
to see that the same complication turns up in the one-dimensional case.
\end{remark}

%
\begin{remark}\label{rkpropMN} The functions $M$ and $N$ involved
in the second and third-order conditions satisfy some usual properties;
see, for example, \citet{dehaanresnick1993}. More specifically, one
can show that there exist nonpositive reals $\rho$ and $\rho'$ such
that $ \alpha$ (resp., $ \beta$) is a regularly varying function of
order $\rho$ (resp., $ \rho'$), that is, $ \alpha(tz)/\alpha(t) \to
z^\rho$
when $t \to\infty$, for each positive $z$.
Besides, $M$ is homogeneous of order $1-\rho$, that is to say $M(r{\mathbf
x})=r^{1-\rho} M({\mathbf x})$, for each positive $r$ and ${\mathbf x}$ with
positive coordinates.
Finally, the function $N$ is homogeneous of order $1-\rho- \rho'$.
\end{remark}

%
\begin{remark}\label{rkAI}
An interesting situation to consider is when the c.d.f. $F$ is in the
domain of attraction of an extreme value distribution $G$ with
independent components, that is, $G = \prod_{j=1}^d G_j$. Such a c.d.f.
is said to have the property of {asymptotic independence}. In this
case, the function $M$ is the limit of the joint tail of the
distribution, and in dimension 2, the coefficient of tail dependence
$\eta$ introduced by \citet{ledfordtawn1996,ledfordtawn1997}
equals $1/(1-\rho)$, where $\rho$ is defined in Remark~\ref{rkpropMN}.
\end{remark}

In this paper, we will handle two sets of assumptions. First consider:
\begin{itemize}[(A2)~~--]
\item[(A2)~~--] the second-order condition is satisfied, so that~(\ref
{eq2ndorder}) holds;
\item[--] the coefficient of regular variation $\rho$ of the function
$\alpha$ defined in~(\ref{eq2ndorder}) is negative;
\item[--] the function $M$ defined in~(\ref{eq2ndorder}) is continuous.
\end{itemize}
These hypotheses allow us to get the asymptotic uniform behavior of
$\hat L_k$, the empirical estimator of $L$ defined by~(\ref
{eq1storderhat}), as detailed in the following proposition.

%
\begin{proposition}\label{propcv-ps-L}
Let ${\mathbf X}_1, \ldots, {\mathbf X}_n$ be independent multivariate random
vectors in $\mathbb{R}^d$ with common joint c.d.f.
$F$ and continuous marginal c.d.f.'s $F_j$ for $j=1,\ldots,d$. Assume that
the set of conditions \textup{(A2)} holds. Suppose further that the
first-order partial derivatives of $L$ (denoted by $ \partial_jL$ for
$j=1,\ldots,d$) exist and that $\partial_j L$ is continuous on the
set of points
$\{{\mathbf x}=(x_1,\dots, x_d) \in\mathbb{R}^d_+\dvtx  x_j >0\}$.
Consider $\hat L_k$ the estimator of $L$ defined by~(\ref
{eq1storderhat}) where $k$ is such that $\sqrt{k} \alpha(n/k) \to
\infty$.
Then as $n$ tends to infinity,
we get
\[
\sup_{0\leq x_1,\ldots, x_d \leq T} \biggl\llvert \frac{1}{\alpha(n/k)} \bigl\{ \hat
L_k ({\mathbf x}) - L({\mathbf x}) \bigr\} -M({\mathbf x}) \biggr\rrvert
\stackrel{\mathbb{P}} {\longrightarrow} 0.
\]
\end{proposition}

Under stronger assumptions, and for some choice of the intermediate
sequence, the asymptotic distribution of the previous stochastic
process can be obtained after multiplication by the rate $\sqrt{k}
\alpha(n/k)$. For a positive $T$, let $D([0,T]^d)$ be the space of
real valued functions that are right-continuous with left-limits.
Now introduce the conditions:
\begin{itemize}[(A3)~~--]
\item[(A3)~~--] the third-order condition is satisfied, so that~(\ref
{eq2ndorder}) and~(\ref{eq3rdorder}) hold;
\item[--] the coefficients of regular variation $\rho$ and $\rho'$ of
the functions $\alpha$ and $\beta$ defined in~(\ref{eq2ndorder})
and~(\ref{eq3rdorder}) are negative;
\item[--] the function $M$ defined in~(\ref{eq2ndorder}) is
differentiable and $N$ defined in~(\ref{eq3rdorder}) is continuous.
\end{itemize}
%

%
\begin{proposition}
\label{propdev-asympt-L}
Assume that the conditions of Proposition~\ref{propcv-ps-L} are
fulfilled and that the set of conditions \textup{(A3)} hold.
Consider $\hat L_k$ the estimator of $L$ defined by~(\ref
{eq1storderhat}) where $k$ is such that $\sqrt{k} \alpha(n/k) \to
\infty$ and $\sqrt{k} \alpha(n/k) \beta(n/k)\to0$.
Then as $n$ tends to infinity,
%
\begin{equation}
\label{eqasympt-dev-L} \sqrt{k} \biggl\{ \hat L_k({\mathbf x}) - L({\mathbf x}) -
\alpha\biggl(\frac
{n}{k}\biggr)M({\mathbf x}) \biggr\} \stackrel{d} {\to}
Z_L({\mathbf x}),
\end{equation}
in $D([0,T]^d)$ for every $T>0$,
where
\begin{equation}
\label{eqZL}Z_L({\mathbf x}):  = W_L({\mathbf x}) - \sum
_{j=1}^d W_L(x_j
{\mathbf e}_j)\partial_jL({\mathbf x}).
\end{equation}
The process $W_L$ above is a continuous centered Gaussian process with
covariance structure
$ {\mathbb E}[W_L({\mathbf x})W_L({\mathbf y})]= \mu\{R({\mathbf x}) \cap R({\mathbf
y})\}$
given in terms of the measure $\mu$ defined by~(\ref{eqmu}) and of
$
R({\mathbf x})=\{{\mathbf u} \in{\mathbb R}^d_+$: there exists $j$ such
that  $0 \leq u_j\leq x_j \}$.
\end{proposition}

%
\begin{remark}\label{rkasymptbias} A difference between the
previous result and Theorem 7.2.2 of \citet{dehaanferreira2006}
consists of the choice of the intermediate sequence that is larger
here. Indeed, we suppose $\llvert \sqrt{k} \alpha(n/k) \rrvert \to\infty$ whereas
they choose $k(n)=o (n^{-2\rho/(1-2\rho)} )$, which
implies $\sqrt{k} \alpha(n/k) \to0$. Our choice requires the more
informative second-order condition~(\ref{eq2ndorder}). A~nonnull
asymptotic bias appears in our framework.
\end{remark}

%
\begin{remark} The conditions on $k$, $\alpha$ and $\beta$
required in Proposition~\ref{propdev-asympt-L} are not too
restrictive: because of the regular variation of $\alpha$ and $\beta$,
they are implied by the choice $k(n)=n^\kappa$, with $
\kappa\in (-\frac{2\rho}{1-2\rho},   -\frac{2(\rho+ \rho
')}{1-2(\rho+ \rho')} )$.
\end{remark}


\section{Bias reduction procedure}\label{secprocedure}
As pointed out in Remark~\ref{rkasymptbias}, a nonnull asymptotic
bias $\alpha({n}/{k})M({\mathbf x})$ appears from Proposition~\ref
{propdev-asympt-L}. The bias reduction procedure will consist of
subtracting the estimated asymptotic bias obtained in Section~\ref
{subsecmethodA}. The key ingredient is the homogeneity of the
functions $L$ and $M$ mentioned in Remarks~\ref{rkhomoL} and \ref
{rkpropMN}. This homogeneity will also provide other constructions to
get rid of the asymptotic bias.

\subsection{Estimation of the asymptotic bias of \texorpdfstring{$\hat L_k$}{hatLk}}
\label{subsecmethodA}
Equation~(\ref{eqasympt-dev-L}) suggests a natural correction of
$\hat L_k$ as soon as an estimator of $\alpha({n}/{k})M({\mathbf x})$ is
available. In order to take advantage of the homogeneity of $L$,
let us introduce a positive scale parameter $a$ which allows to
contract or to dilate the observed points. We denote
%
\begin{equation}
\label{eqLa-hat} \hat L_{k,a}({\mathbf x}):=a^{-1}\hat
L_k(a{\mathbf x})
\end{equation}
and
%
\begin{equation}
\label{eqDelta} \hat \Delta_{k,a}({\mathbf x}):=\hat L_{k,a}({\mathbf
x})-\hat L_{k}({\mathbf x}).
\end{equation}

From~(\ref{eqasympt-dev-L}) one gets
%
\begin{equation}
\label{eqdev-asympt-La-L} \sqrt{k} \biggl\{\hat L_{k,a}({\mathbf x}) -L({\mathbf x}) -
\alpha\biggl(\frac
{n}{k}\biggr) a^{-\rho} M({\mathbf x}) \biggr\}
\stackrel{d} {\to} a^{-1} Z_L(a{\mathbf x}),
\end{equation}
in $D([0,T]^d)$ for every $T>0$. Equations~(\ref{eqDelta}) and
Proposition~\ref{propcv-ps-L} yield as $n$ tends to infinity,
%
\begin{equation}
\label{eqDeltaCV} \frac{ \hat\Delta_{k,a}({\mathbf x})}{\alpha(n/k)} \overset {{\mathbb P}} {\longrightarrow}
\bigl(a^{-\rho} -1\bigr)M({\mathbf x}).
\end{equation}
Fixing $a$ such that $a^{-\rho} -1=1$, a natural estimator of the
asymptotic bias of $\hat L_k({\mathbf x})$ is thus
$ \hat\Delta_{k, 2^{-1/{ \hat\rho}}} ({\mathbf x})$,
where $\hat\rho$ is an estimator of $\rho$. Recall that the unknown
parameter $\rho$ is the regular variation index of the function
$\alpha$ involved in the-second order condition.
Let $k_\rho$ be an intermediate sequence that represents the number of
order statistics used in the estimator~$\hat\rho$. Assume that
$k_\rho\gg k$ where $k=k(n)$ is the sequence used in Proposition~\ref
{propdev-asympt-L}. A first asymptotically unbiased estimator of
$L({\mathbf x})$ can be defined as
%
\begin{equation}
\label{defestimA}
\mathring L_{k, 1, k_\rho}({\mathbf x}):= \hat L_k({\mathbf
x}) - \hat\Delta _{k, 2^{-1/{\hat\rho}}} ({\mathbf x}).
\end{equation}
The asymptotic behavior of this estimator is provided in Theorem~\ref{thmbiasB} and Remark~\ref{rkCP}. We refer the reader to
Section~\ref{secrho} for more details concerning the estimation of
$\rho$.

\subsection{Estimation of the asymptotic bias of \texorpdfstring{$\hat L_{k,a}$}{hatLk,a}}\label{subsecmethodB}
The previous construction can be easily generalized by correcting the
estimator $\hat L_{k,a}$ instead of $\hat L_{k}$. Indeed, from~(\ref
{eqdev-asympt-La-L}) one can see that the asymptotic bias of $\hat
L_{k,a}({\mathbf x})$ is $ \alpha(\frac{n}{k}) a^{-\rho} M({\mathbf x})$.
Recall that
when $n$ tends to infinity, one has for any positive real $b$,
\[
\frac{ \hat\Delta_{k,b}({\mathbf x})}{\alpha(n/k)} \overset {{\mathbb P}} {\longrightarrow}
\bigl(b^{-\rho} -1\bigr)M({\mathbf x}). %
\]
Thus fixing $b$ such that $b^{-\rho} -1=a^{-\rho}$ will help to
cancel the asymptotic bias. It yields the following asymptotically
unbiased estimator of $L$:
%
\begin{equation}
\label{defestimB} \mathring L_{k, a, k_\rho}({\mathbf x}):= \hat L_{k,a}({\mathbf
x}) - \hat \Delta_{k, (a^{-\hat\rho} +1)^{-1/{\hat\rho}}} ({\mathbf x}).
\end{equation}
%

\begin{theorem}
\label{thmbiasB}
Assume that the conditions of Proposition~\ref{propdev-asympt-L} are
fulfilled, and consider the estimator of $L$ defined by~(\ref{defestimB}).
Let $k_\rho$ be an intermediate sequence such that $\sqrt{k_\rho
}\alpha(n/k_\rho)(\hat\rho-\rho)$ converges in distribution.
Suppose also that $k$ is such that $k=o(k_\rho)$, $\sqrt{k} \alpha
(n/k) \to\infty$ and $\sqrt{k} \alpha(n/k) \beta(n/k) \to 0$.
Under these assumptions, as $n$ tends to infinity,
%
\begin{equation}
\label{eqdev-asympt-B} \sqrt{k} \bigl\{ \mathring L_{k, a, k_\rho}({\mathbf x}) -L({\mathbf x})
\bigr\} \stackrel{d} {\to} \mathring Y_{a}({\mathbf x}),
\end{equation}
in $D([0,T]^d)$ for every $T>0$, where $ \mathring Y_{a}$ is a
continuous centered Gaussian process defined by
\[
\mathring Y_{a} ({\mathbf x}):= Z_L({\mathbf x}) -
b^{-1} Z_L( b {\mathbf x} ) + a^{-1}
Z_L(a {\mathbf x}) %
\]
with covariance
$
{\mathbb E}[ \mathring Y_a({\mathbf x}) \mathring Y_a({\mathbf y})] = {\mathbb
E}[Z_L({\mathbf x}) Z_L({\mathbf y})]  (1 - b^{-1/2} + a^{-1/2} )^2$
and $b = (a^{-\rho} +1)^{-1/\rho}$.
\end{theorem}

%
\begin{remark}\label{rkhypo4}
The assumption that $\sqrt{k_\rho}\alpha(n/k_\rho)(\hat\rho-\rho
)$ converges in distribution will be reconsidered in Section~\ref{secrho}.
\end{remark}

%
\begin{remark}\label{rkCP}
Theorem~\ref{thmbiasB} remains true when $a=1$ and thus characterizes
the asymptotic behavior of the estimator given in~(\ref{defestimA}).
For this particular choice of $a$, the covariance reduces to $ {\mathbb
E}[Z_L({\mathbf x}) Z_L({\mathbf y})](2-2^{1/{2\rho}})^2$.
\end{remark}
%
\subsection{An alternative estimation of the asymptotic bias of \texorpdfstring{$\hat L_{k,a}$}{hatLk,a}}\label{subsecmethodC}
The procedure of bias reduction introduced in the previous section
requires the estimation of the second-order parameter~$\rho$. It is
actually possible to avoid it, making use of combinations of estimators
of $L$. The asymptotic bias of $\hat L_{k,a}({\mathbf x})$ is $ \alpha
(\frac{n}{k}) a^{-\rho} M({\mathbf x})$, as already noted from~(\ref
{eqdev-asympt-La-L}).
Making use of~(\ref{eqDeltaCV}) and homogeneity of $M$, one gets as
$n$ tends to infinity,
\[
\frac{ \hat\Delta_{k_\rho,a}(a {\mathbf x})} { \hat\Delta_{k_\rho,a}(a {\mathbf x}) -a \hat\Delta_{k_\rho,a}({\mathbf x})} \stackrel{ \mathbb{P} } {\longrightarrow} \frac{a^{-\rho}}{a^{-\rho}-1}, %
\]
for any intermediate sequence $k_\rho$ that satisfies $\sqrt{k_\rho
}\alpha(n/k_\rho)\to\infty$.
The expression
\[
\hat\Delta_{k,a}({\mathbf x}) \frac{ \hat\Delta_{k_\rho,a}(a {\mathbf x})
}{ \hat\Delta_{k_\rho,a}(a {\mathbf x}) -a \hat\Delta_{k_\rho,a}({\mathbf
x}) } %
\]
can thus be used as an estimator of the asymptotic bias of $\hat
L_{k,a}({\mathbf x})$. After simplifications,
this leads to a new family of asymptotically\vspace*{1pt} unbiased estimators of
$L({\mathbf x})$ by substracting the estimated bias from $\hat L_{k,a}({\mathbf
x})$, namely
%
\begin{equation}
\label{defestimC} \tilde L_{k, a, k_\rho}({\mathbf x}) =\frac{\hat L_k({\mathbf x}) \hat\Delta
_{k_\rho,a}(a {\mathbf x})-\hat L_k(a {\mathbf x}) \hat\Delta_{k_\rho,a}({\mathbf x})}{\hat\Delta_{k_\rho,a}(a {\mathbf x}) -a \hat\Delta
_{k_\rho,a}({\mathbf x}) },
\end{equation}
which is well defined for any real number $a$ such that $0<a< 1$.

%
\begin{theorem}
\label{thmbias-Lu}
Assume that the conditions of Proposition~\ref{propdev-asympt-L} are
fulfilled, and consider the estimator of $L$ defined by~(\ref{defestimC}).
Let $k_\rho$ be an intermediate sequence such that $\sqrt{k_\rho
}\alpha(n/k_\rho)(\hat\rho-\rho)$ converges in distribution.
Suppose also that $k$ is such that $k=o(k_\rho)$, $\sqrt{k} \alpha
(n/k) \to\infty$, $\sqrt{k}=O(\sqrt{k_\rho}\alpha(n/k_\rho))$
and $\sqrt{k} \alpha(n/k) \beta(n/k) \to 0$. Assume moreover that
the function $M$ never vanishes except on the axes.
Then, as $n$ tends to infinity,
%
\begin{equation}
\label{eqdev-asympt-unbiased} \sqrt{k} \bigl\{ \tilde L_{k, a, k_\rho}({\mathbf x}) -L({\mathbf x})
\bigr\} \stackrel{d} {\to} \tilde Y_a({\mathbf x}), 
\end{equation}
in $D([\varepsilon,T]^d)$ for every $\varepsilon>0$ and $T>0$, where $
\tilde Y_{a}$ is a continuous centered Gaussian process with covariance
$
{\mathbb E}[ \tilde Y_a({\mathbf x}) \tilde Y_a({\mathbf y})]$ given by
$ {\mathbb E}[Z_L({\mathbf x}) Z_L({\mathbf y})]\times\break (a^{-\rho} -1)^{-2}  (
a^{-\rho} - a^{-1/2} )^2$.
\end{theorem}

%
\begin{remark}
The covariance function specified above is decreasing with respect to
the parameter~$a$ for any fixed value of $\rho$. This suggests at
first glance to choose $a$ close to 1 in order to reduce the asymptotic
variance of $ \tilde Y_{a}$, but this would give a degenerate form of
(\ref{defestimC}). See Section~\ref{secsimu} for practical
considerations for the choice of $a$.
\end{remark}


\section{Theoretical examples}\label{secexamples}
The aim of this section is to furnish several multivariate
distributions that satisfy the third-order condition~(\ref
{eq3rdorder}). For the sake of simplicity, expressions are displayed
in the bivariate setting.
We start by focusing on heavy-tailed margins. In this case, a first
possible step to get the pointwise convergence is to obtain, for well-chosen
positive reals $p$ and $q$, an expansion (for $t$ tending to infinity)
of the form
\begin{eqnarray*}
&& t {\mathbb P}\bigl(X>t^px \mbox{ or }  Y>t^q y\bigr)
\\
&&\qquad  =
T_1(x,y) + \alpha(t) T_2(x,y) + \alpha(t) \beta(t)
T_3(x,y) + o\bigl(\alpha(t) \beta(t) \bigr), %
\end{eqnarray*}
with $T_1(1,1)>0$. One can then identify each term involved in~(\ref
{eq3rdorder}) as follows:
\begin{eqnarray*}
L(x,y) &=& T_1\bigl(a(x), b(y)\bigr),\qquad M(x,y) = T_2
\bigl(a(x), b(y)\bigr)\quad\mbox{and}\quad
\nonumber\\[-8pt]\\[-8pt]\nonumber
N(x,y) &=&  T_3\bigl(a(x), b(y) \bigr),
\end{eqnarray*}
where
\[
a(x) = x^{-p} \bigl\{T_1(1,+\infty) \bigr\}^p, \qquad  b(x) = x^{-q} \bigl\{T_1(+\infty,1) \bigr
\}^q. %
\]
Applying \citeauthor{resnick1986} [(\citeyear{resnick1986}), Corollary 5.18], one can check that
in such a framework a form of the bivariate extreme value distribution
$G$ is given by
\[
G(x,y)=\exp \biggl(-\frac{T_1(x,y)}{T_1(1,1)} \biggr).
\]

\subsection{Powered norm densities}\label{subsecResnick-style}
Following the idea of \citeauthor{resnick1986} [(\citeyear{resnick1986}), pages 276~and~286]
consider first a norm $\llVert   \cdot\rrVert  $, and a cone $
\mathcal{D}$ of $\mathbb{R}^2$, that is to say, a~set such that if $(x,y)\in
\mathcal{D}$, then $(tx,ty)\in\mathcal{D}$ for every positive~$t$.
Without loss of generality, suppose that $(1,1)\in\mathcal{D}$. Let
$(X,Y)$ be a bivariate random vector with probability density function
given by
\[
f(x,y):=\frac{c {\mathbf1}_{\mathcal{D}}(x,y)}{(1+\llVert  (x,y)^T\rrVert  ^\alpha
)^{\beta}},
\]
where $c$ is a normalizing positive constant and where $\alpha$ and
$\beta$ are some positive real numbers such that $\alpha\beta>2$.
Set $A_\mathcal{D}(x,y):=\{(u,v)\in\mathcal{D}\dvtx  u>x$ or
$v>y\}$, and define $p:=(\alpha\beta-2)^{-1}$.
One can check that for $j=1,2,3$,
\begin{eqnarray*}
T_j(x,y)&=&\iint_{A_\mathcal{D}(x,y)}\frac{c   c_j  \,du \,dv}{\llVert
(u,v)^T\rrVert  ^{\alpha(\beta+j-1)}},
\end{eqnarray*}
where $c_1=1$, $c_2=-\beta$ and $c_3=\beta(\beta+1)/2$.
The functions $M$ and $N$ are homogeneous with order given through
$\rho=\rho^\prime=-\alpha p$.

Let us discuss some particular choices of the norm:
\begin{itemize}[--]
\item[--] For the $L^1$-norm and $\alpha=1$, the model coincides with
the bivariate Pareto of type II distribution,\vspace*{1pt} denoted by $\operatorname{BPII}(\beta)$
in this paper, and referred to as MP$^{(2)}(\mathit{II})(0,1, \beta-2)$
in \citet{KBJ2000}, page~604. In this case, $p=q=( \beta- 2)^{-1}$, and
$
L(x,y)= x+y -(x^{-p} + y^{-p})^{-1/p}$.
The latter s.t.d.f. is known as the negative logistic model, introduced by
\citet{joe1990}; see also \citet{beirlantgoegebeursegersteugels2004}, page 307.
\item[--] When the Euclidean norm is chosen, one recovers the bivariate
Cauchy distribution for $\alpha=2$, $\beta=3/2$ and $p=1$. On the
positive quadrant, that means for $\mathcal{D}=\mathbb{R}_+^2$, we
have $c=2/\pi$, $T_1(u,v)=c(u^{-2}+v^{-2})^{1/2}$ and $a(x)=b(x)=c/x$.
On the whole plane, which means that $\mathcal{D}=\mathbb{R}^2$, we
get $c=1/(2\pi)$, $T_1(u,v)=c \{u^{-1} + v^{-1}+
(u^{-2}+v^{-2})^{1/2} \}$ and $a(x)=b(x)=2c/x$. This can also be
seen as a particular case of the following item.
\item[--] The Student distributions with Pearson correlation
coefficient $\theta$ arise choosing the norm $\llVert  (x,y)^T\rrVert  =\nu
^{-1/2}(x^2-2\theta xy+y^2)^{1/2}$, for a positive real number $\nu$,
$\alpha=2$, $\beta=(\nu+2)/2$ and $p=\nu^{-1}$. In this case, the
integral form of the function $T_1$ cannot be totally simplified, and
one classically writes the s.t.d.f. as
\begin{eqnarray*}
L(x,y) &=& (x+y) \biggl[ \frac{y}{x+y}F_{\nu+1} \biggl\{
\frac{(y/x)^{1/\nu}-\theta}{\sqrt{1-\theta^2}} \sqrt{\nu+1} \biggr\}
\\
&&\hspace*{38pt}{}  + \frac{x}{x+y}F_{\nu+1}
\biggl\{\frac{(x/y)^{1/\nu}-\theta}{\sqrt
{1-\theta^2}} \sqrt{\nu+1} \biggr\} \biggr], %
\end{eqnarray*}
where $F_{\nu+1}$ is the c.d.f. of the univariate Student
distribution with $\nu+1$ degrees of freedom. This dependence
structure is also obtained for some elliptical models; see, for example,
[\citet{krajina2012}, page 1813] and next subsection.
\item[--] Other choices for the norm would lead to other distributions.
Note that one can also relax the symmetry condition, considering, for
instance, the Mahalanobis pseudo-norm defined by $\llVert  (x,y)^T\rrVert
^2=(x/\sigma)^2- 2\rho(x/\sigma)(y/\tau)+(y/\tau)^2$ for a real
number $\rho$ such that $\llvert \rho\rrvert <1$ and some positive real numbers
$\sigma$ and $\tau$.
\end{itemize}

\subsection{Elliptical distributions}\label{subsecelliptical}
Consider the usual representation of the centered elliptical
distribution $(X,Y)^T=R{\mathbf AU}$, in terms of a positive random
variable $R$, a $2\times2$ matrix $\mathbf A$ such that ${\bolds\Sigma
}={\mathbf A A}^T$ is of full rank, and a bivariate random vector $\mathbf U$
independent of $R$, uniformly distributed on the unit circle of the plane.
Assume that $R$ has a probability density function denoted by $g_R$.
One can then express the probability density function of $(X,Y)$ as
\[
f(x,y):=\frac{1}{\llvert \operatorname{det}\mathbf A\rrvert }g_R \bigl\{ (x, y) {\bolds\Sigma
}^{-1} (x,y)^T \bigr\}. %
\]
A sufficient condition to satisfy~(\ref{eq3rdorder}) is to assume
that the distribution of $R$ belongs to the Hall and Welsch class
[\citet{hallwelsh1985}], namely,
\[
{\mathbb P}(R>r) = c r^{-1/\gamma} \bigl\{ 1 + D_1
r^{\rho/\gamma} + D_2 r^{(\rho+ \rho_1)/\gamma} + o\bigl(r^{(\rho+ \rho_1)/\gamma
}
\bigr) \bigr\}, %
\]
with positive real $c$, nonnull reals $D_1$ and $ D_2$ and negative
reals $\rho$ and $ \rho_1$.

One can check that, for $j=1,2,3$,
\begin{eqnarray*}
T_j(x,y)&=&\frac{c}{2 \pi\gamma\llvert \operatorname{det}\mathbf A\rrvert } \iint_{\{(u,v)\dvtx  u
> x~\mathrm{or}~v >y\}}\frac{du \,dv}{\{(u,v) {\bolds\Sigma}^{-1}(u,v)^T\}
^{1+1/(2\gamma) + p_j}},
\end{eqnarray*}
where $p_1=0, p_2 = - \rho/(2\gamma)$ and $p_3 = - (\rho+ \rho
_1)/(2\gamma)$.

Assuming for simplicity that ${\bolds\Sigma} = {1\ \  \theta\choose \theta\ \ 1}$, the s.t.d.f. can be written as
\begin{eqnarray*}
L(x,y) &=& (x+y) \biggl[ \frac{y}{x+y}F_{1/\gamma+1} \biggl\{
\frac
{(y/x)^{\gamma}-\theta}{\sqrt{1-\theta^2}} \sqrt{1/\gamma+1} \biggr\}
\\
&&\hspace*{37pt}{} + \frac{x}{x+y}F_{1/\gamma+1}
\biggl\{\frac{(x/y)^{\gamma
}-\theta}{\sqrt{1-\theta^2}} \sqrt{1/\gamma+1} \biggr\} \biggr],
\end{eqnarray*}
which is the form already obtained for the Student distribution in
Section~\ref{subsecResnick-style} for $\nu= 1/\gamma$. See \citet
{demartamcneil2005} for more details. Note finally that for a general
matrix $\bolds\Sigma$ and the special case $g_R(r)= c(1+r^\alpha
)^{-\beta}$, one recovers the Mahalanobis pseudo-norm already
mentioned in the previous subsection.

When dealing with margins that are \textit{not} heavy tailed, the calculus
is done directly from~(\ref{eq2ndorder}).
The last two examples of bivariate distributions have short and light
tailed margins, respectively.

\subsection{Archimax distributions}\label{subsecArchimax}
Consider the bivariate c.d.f. defined for each $0 \leq u,v \leq1$ by
%
\begin{equation}
\label{eqArchimax} F(u,v) = \bigl\{1 + L\bigl(u^{-1}-1,v^{-1} -1
\bigr) \bigr\}^{-1},
\end{equation}
given in terms of a s.t.d.f. $L$.
This distribution has standard uniform univariate margins and
corresponds to a particular case of Archimax bivariate copulas
introduced in \citet{caperaafougeresgenest2000}, in which the
function $\phi(t) = t^{-1} -1$ is the Clayton Archimedean generator
with index 1. Expanding the left-hand side term of~(\ref{eq2ndorder})
leads to, as $t$ tends to infinity,
\[
t \bigl\{ 1- F \bigl(1-t^{-1} x, 1- t^{-1} y \bigr) \bigr\}
= L(x,y) +t^{-1} M(x,y) + t^{-2} N(x,y) + o
\bigl(t^{-2} \bigr), %
\]
where
\begin{eqnarray*}
M(x,y) &:=& x^2 \partial_1 L(x,y)+y^2
\partial_2 L(x,y) -L^2(x,y),
\\
N(x,y) &:=& x^4/2 \partial^2_{11}L(x,y)
+x^2y^2 \partial^2_{12}L(x,y)+y^4/2
\partial^2_{22}L(x,y)
\\
&&{}+ L^3(x,y) + \bigl(x^3- 2
x^2 L(x,y) \bigr) \partial_1 L(x,y)
\\
&&{} + \bigl(y^3 - 2 y^2 L(x,y) \bigr) \partial_2
L(x,y).
\end{eqnarray*}
This allows us to identify $\rho= \rho^\prime= -1$. Above, the
notation $\partial_{ij}L$ stands for
$\partial^2 L / (\partial x_i\, \partial x_j)$.

\subsection{Multivariate symmetric logistic distributions}\label{subsecasymlog}
Consider the c.d.f. defined by
%
\begin{equation}
\label{eqasymlog} F(x,y) = \exp \bigl\{ - \bigl(e^{-x/s} + e^{-y/s}
\bigr)^s \bigr\},
\end{equation}
for each $x,y \in{\mathbb R}$, which corresponds to the bivariate
extreme value distribution with Gumbel univariate\vspace*{1pt} margins
$F_1(x)=F_2(x)=\exp\{ - e^{-x} \} $ and symmetric logistic s.t.d.f.
$L(x,y)= (x^{1/s} + y^{1/s})^s$,
where $0 < s \leq1$. This distribution was introduced in \citet
{tawn1988}; see, for example, \citet{beirlantgoegebeursegersteugels2004}, page~304.
Expanding $
t    [ 1- F \{F_1^{-1}(1-t^{-1} x),   F_2^{-1}(1- t^{-1} y)
 \} ]$
leads to 
\[
L(x,y) +t^{-1} M(x,y) + t^{-2} N(x,y) + o
\bigl(t^{-2} \bigr), %
\]
%
where
\begin{eqnarray*}
M(x,y) &:=& \tfrac{1}{2} \bigl(x x^{1/s}+y y^{1/s}\bigr)
\bigl\{L(x,y)\bigr\}^{1-1/s}- \tfrac
{1}{2}\bigl\{L(x,y)\bigr
\}^2,
\\
N(x,y)&:=&\frac{1}{3} \bigl(x^2 x^{1/s}+y^2
y^{1/s}\bigr)\bigl\{L(x,y)\bigr\}^{1-1/s}
\\
&&{} +\frac{1-s}{8s}
(xy)^{1/s}(x -y)^2 \bigl\{L(x,y)\bigr\}^{1-2/s}
\\
&&{} +\frac{1}{3! }\bigl\{L(x,y)\bigr\}^3 -
\frac{1}{2} \bigl(x x^{1/s}+y y^{1/s}\bigr)\bigl\{L(x,y)
\bigr\}^{2-1/s}.
\end{eqnarray*}
This allows us to identify $\rho= \rho^\prime= -1$.
The identification of second and third-order terms has previously be
derived by \citet{ledfordtawn1997}.

\section{Finite sample performances}\label{secsimu}
The purpose of this section is to evaluate the performance of the
estimators of $L$ introduced in Section~\ref{secprocedure}. For
simplicity, we will focus on dimension 2, and simulate samples
from the distributions presented in Section~\ref{secexamples}.
Thanks to the homogeneity property, one can focus on the estimation of
$t \mapsto L(1-t,t)$ for $0 \leq t \leq1$, which coincides with the
Pickands dependence function $A$; see, for example, \citet
{beirlantgoegebeursegersteugels2004}, page 267. Considering first
the estimation at $t=1/2$ leads to the definition of aggregated
versions of our estimators. These new estimators will be both compared
in terms of $L^1$-errors for $L$ or associated level curves.

\subsection{Estimators in practice}\label{subseccorrections}
\begin{figure}

\includegraphics{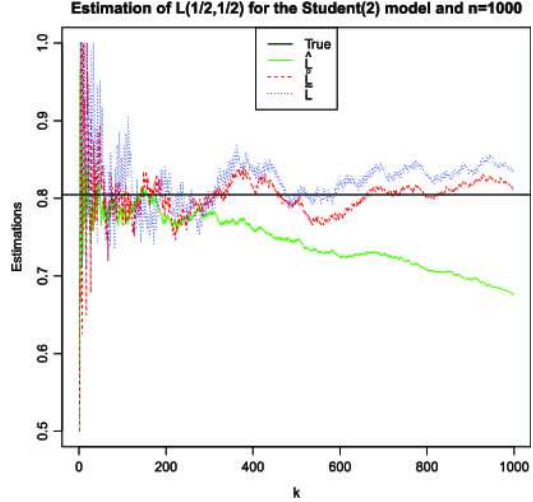}

\caption{Estimation of $L(1/2,1/2)$ for the bivariate $\operatorname{Student}(2)$ law
based on a sample of size 1000.}
\label{sec-simu1-stu2}
\end{figure}

Let us start with the estimation of $L(1/2,1/2)$ for the bivariate
Student distribution with 2 degrees of freedom. This model is a
particular case of Sections~\ref{subsecResnick-style} and \ref
{subsecelliptical}. For one sample of size 1000, Figure~\ref
{sec-simu1-stu2} gives, as functions of $k$, the estimation of $L$ at
point $(1/2,1/2)$ by $\hat L_{k}$, $\mathring L_{k}$ and $\tilde L_{k}$,
respectively, defined by~(\ref{eq1storderhat}),~(\ref{defestimB})
and~(\ref{defestimC}). For the last two estimators, the parameters
have been tuned as follows:
$a=0.4$, $k_\rho= 990$ and $\rho$ estimated using~(\ref
{defrho-hat}) with $a=r=0.4$. These values have been empirically
selected based on intensive simulation, and will be kept throughout the paper.
%
One can check from Figure~\ref{sec-simu1-stu2} that the empirical
estimator $\hat L_{k}$ behaves fairly well in terms of bias for small
values of~$k$. Besides, the bias is efficiently corrected by the two
estimators $\mathring L_{k}$ and $\tilde L_{k}$. Since the bias almost
vanishes along the range of $k$,
one can think about reducing the variance through an aggregation in $k$
(via mean or median) of $\mathring L_{k}$ or $\tilde L_{k}$. This leads
us to consider
the two following estimators:
\begin{eqnarray*}
\mathring L_{\mathrm{agg}} &:=& \operatorname{Median}(\mathring L_{k}, k=1,
\ldots, \kappa_n),
\\
\tilde L_{\mathrm{agg}}&:=&\operatorname{Median}(\tilde L_{k}, k=1, \ldots,
\kappa _n), %
\end{eqnarray*}
where\vspace*{1pt} $n$ is the sample size and $\kappa_n$ is an appropriate fraction
of $n$. Their performance will be compared to those of the family $\{
\hat L_k, k=1,\ldots,n-1\}$. Simplified notation $\{\hat L_k, k\}$
will be used instead of $\{\hat L_k, k=1,\ldots,n-1\}$. Because any
s.t.d.f.~$L$ satisfies $\max(t,1-t) \leq L(1-t,t) \leq1$, the competitors
have been corrected so that they satisfy the same inequalities.

%
\begin{remark} If $\kappa_n$ satisfies the condition imposed on
$k_n$ in Theorems~\ref{thmbiasB} and \ref{thmbias-Lu}, then the
aggregated estimators $\mathring L_{\mathrm{agg}}$ and $\tilde L_{\mathrm{agg}}$ would inherit the asymptotic properties of $\mathring L_{k}$
and $\tilde L_{k}$. Indeed, all the estimators jointly converge, since
they are based on a single process.
\end{remark}

%
\begin{remark}\label{remmixture} In the following simulation
study, $\kappa_n$ is arbitrarily fixed to $n-1$. Such a choice is open
to criticism since it does not satisfy the theoretical assumptions
mentioned in the previous remark. But it is motivated here by the fact
that the bias happened to be efficiently corrected, even for very large
values of $k$, as already illustrated on Figure~\ref{sec-simu1-stu2}.
Note, however, that such a choice would not be systematically the right
one. In presence of more complex models such as mixtures, $\kappa_n$
should not exceed the size of the subpopulation with heaviest tail. To
illustrate this point, take, for example, the bivariate c.d.f. $F=pG +
(1-p)H$, where
$G$ is the c.d.f. of the bivariate $\operatorname{BPII}(3)$ model, and $H$ is the uniform
c.d.f. on $[0,1]^2$. Then the s.t.d.f. is $L(x,y)=x+y-(1/x + 1/y)^{-1}$, and
only $p\%$ of the data belong to the targeted domain of attraction, so
$\kappa_n$ should not exceed $pn$.
\end{remark}

\begin{figure}[b]

\includegraphics{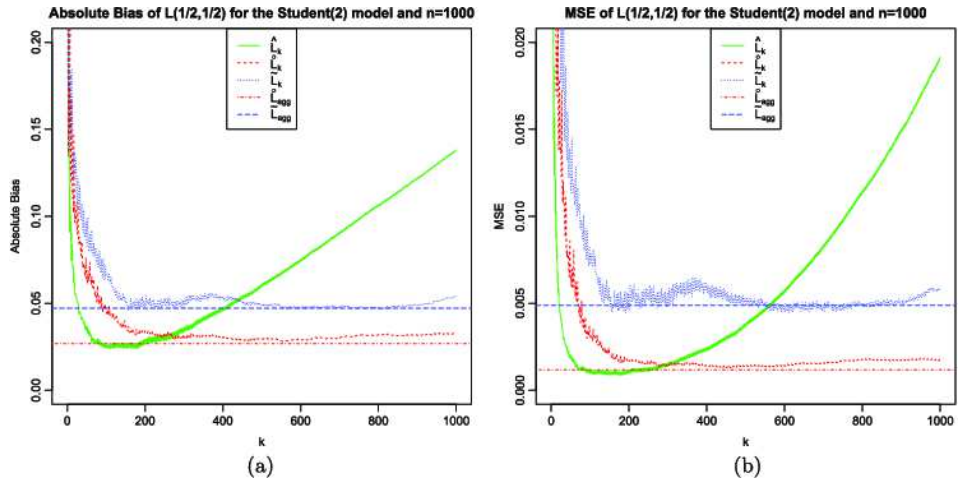}

\caption{\textup{(a)} ABias, \textup{(b)} MSE for the estimation of $L(1/2,1/2)$ in the
bivariate $\operatorname{Student}(2)$ model when $n=1000$ as a function of $k$.}\label{figabias-mse-st}
\end{figure}

Classical criteria of quality of an estimator $\hat\theta$ of $\theta
$ are the absolute bias (ABias) and the mean square error (MSE) defined by
\begin{eqnarray*}
\operatorname{ABias}&=&\frac{1}{N}\sum_{i=1}^N
\bigl\llvert \hat{\theta}^{(i)}-\theta\bigr\rrvert,
\\
\operatorname{MSE}&=&\frac{1}{N}\sum_{i=1}^N
\bigl(\hat{\theta}^{(i)}-\theta\bigr)^2,
\end{eqnarray*}
where $N$ is the number of replicates of the experiment and $\hat
{\theta}^{(i)}$ is the estimate from the $i$th sample.
Note that what we call \textit{Abias} is also referred as \textit{MAE} (for
Mean Absolute Error) in the literature.
Figure~\ref{figabias-mse-st} plots these criteria in the estimation
of $L(1/2,1/2)$ for the bivariate $\operatorname{Student}(2)$ model when $n=1000$ and $N=200$.
%
Figure~\ref{figabias-mse-st} exhibits the strong dependence of the
behavior of $\hat{L}_k$ in terms of $k$, as well as the efficiency of
the bias correction procedures. The estimator $\mathring L_k$ given
by~(\ref{defestimB}) outperforms the estimator $\tilde L_k$ defined
by~(\ref{defestimC}), no matter the value of $k$. Moreover, the ABias
and MSE curves associated to $\mathring L_k$ almost reach the minimum
of those of~$\hat{L}_k$. Finally, the aggregated version $\mathring
L_{\mathrm{agg}}$ answers surprisingly well to the estimation problem of
the s.t.d.f. $L$. First, its performance is similar to the best reachable
from the original estimator $\hat L_k$. Second, it gets rid of the
delicate choice of a threshold $k$ (or would at least simplify this
choice; see Remark \ref{remmixture}).
These comparisons have also been made for five other models obtained
from Section~\ref{secexamples}. The results are very similar to the
ones obtained for the bivariate $\operatorname{Student}(2)$ distribution and are
therefore not presented.

\subsection{Comparisons in terms of $L^1$-error for $L$}\label{subsecL1}

\begin{figure}

\includegraphics{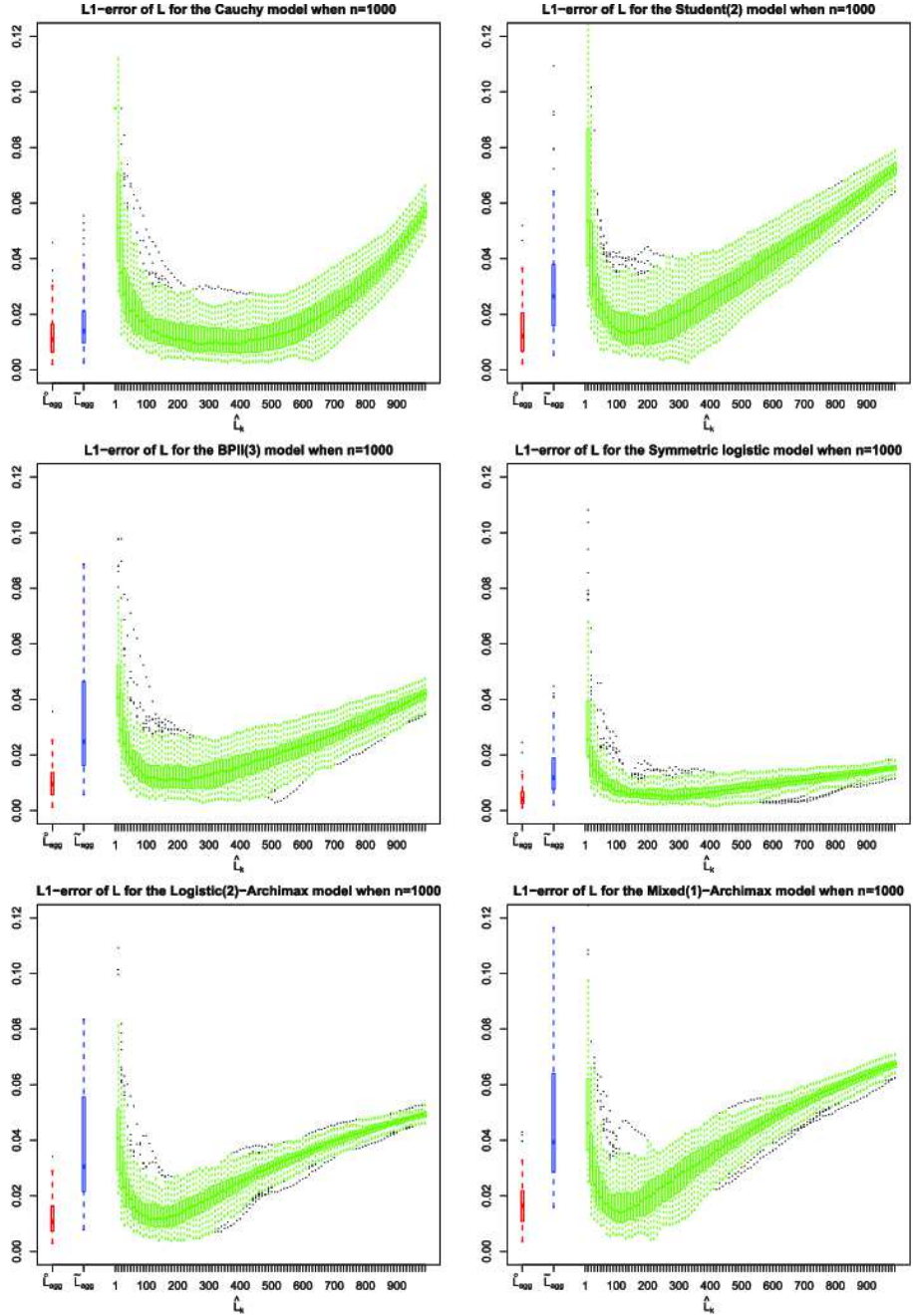}

\caption{Boxplot of the $L^1$-error of function $L$ for the estimators
$\protect\mathring L_{\mathrm{agg}}, \tilde L_{\mathrm{agg}}$ and $\{\hat L_k, k\}$.
First row: bivariate Cauchy model (left) and bivariate $\operatorname{Student}(2)$ model (right).
Second row: bivariate $\operatorname{BPII}(3)$ model (left) and bivariate Symmetric
logistic model (right).
Third row: bivariate Archimax model with logistic (left) and mixed
generator (right).}\label{figl1-norm-L}
\end{figure}

The comparisons are now handled not only at a single point, but for the
whole function using an $L^1$-error defined as follows:
%
\begin{equation}
\label{eqnorm1} \frac{1}{T+1}\sum_{t=1}^T
\biggl\llvert \hat L \biggl(1-\frac{t}{T},\frac
{t}{T} \biggr)- L
\biggl(1-\frac{t}{T},\frac{t}{T} \biggr)\biggr\rrvert,
\end{equation}
where $T$ is the size of the subdivision of $[0,1]$. Figure~\ref{figl1-norm-L} gives the boxplots based on $N=100$ realizations of
$\mathring L_{\mathrm{agg}}, \tilde L_{\mathrm{agg}}$ and $\{\hat L_k, k\}$
for $T=30$ in the case of six bivariate models:
\begin{itemize}
\item First row: Cauchy and $\operatorname{Student}(2)$ models;
\item Second row: $\operatorname{BPII}(3)$ model and Symmetric logistic model with $s=1/3$;
\item Third row: Archimax model with logistic generator
$L(x,y)=(x^2+y^2)^{1/2}$ and mixed generator $L(x,y)=(x^2+y^2+xy)/(x+y)$.
\end{itemize}

%

As already observed in Figure~\ref{figabias-mse-st}, the estimator
$\mathring L_{\mathrm{agg}}$ is again very competitive compared to the
best element of $\{\hat L_k, k\}$, no matter the choice of model.
Recall that the value of $k$ leading to the best $\hat L_k$ depends
crucially on the model and is consequently unknown in practice, which
invites any users to apply this new procedure.

The estimator $ \tilde L_{\mathrm{agg}}$ is definitely less competitive
compared to $\mathring L_{\mathrm{agg}}$. Given these results we will not
pursue with the $ \tilde L_{\mathrm{agg}}$ estimator in the rest of this
paper, and will focus our attention on the behavior of $\mathring
L_{\mathrm{agg}}$.

\subsection{Comparisons between \texorpdfstring{$\protect\mathring{L}_{\mathrm{agg}}$}{mathringLagg}, 
a convex version of \texorpdfstring{$\protect\mathring{L}_{\mathrm{agg}}$}{mathringLagg}, and Peng's estimator} \label{subseccompar}
A natural step is now to compare the performance of our best estimator
$\mathring L_{\mathrm{agg}}$ with
an existing competitor, recently introduced by \citet{peng2010}. In
his work, Peng provides a data-driven method
which chooses the threshold via estimating a s.t.d.f. Another interesting
task is to compare
$\mathring L_{\mathrm{agg}}$ with a convexified version of itself, since
any s.t.d.f. is a convex function; see, for example, \citeauthor{beirlantgoegebeursegersteugels2004}
[(\citeyear{beirlantgoegebeursegersteugels2004}), Section~8.2.2] or \citeauthor{dehaanferreira2006} [(\citeyear{dehaanferreira2006}), Section~6.1.5].
Note that a general convexification procedure has been proposed in
dimension 2 by \citet{filsguillousegers2008}; see also some
alternative suggestions in \citet{bucherdettevolgushev2011}.

In order to take maximal advantage from this simulation study, the
three different models implemented have been considered in two versions
for each: the first model is the Gaussian one, simulated with Pearson's
correlation coefficient $ \pm0.5$. The Gaussian model is a particular
case of elliptical distributions (see Section~\ref
{subsecelliptical}), which illustrates the asymptotic independent
situation; cf. Remark~\ref{rkAI}.
The second model is the bivariate Symmetric logistic one, introduced in
Section~\ref{subsecasymlog}, with two different strengths of
dependence (close to independence on the left column, stronger
dependence on the right column). The third model is the bivariate
Student family, introduced in Sections~\ref{subsecResnick-style}
and~\ref{subsecelliptical} as a particular case. Two strengths of
dependence have also been chosen, close to asymptotic independence on
the left column and stronger dependence on the right column.

\begin{figure}

\includegraphics{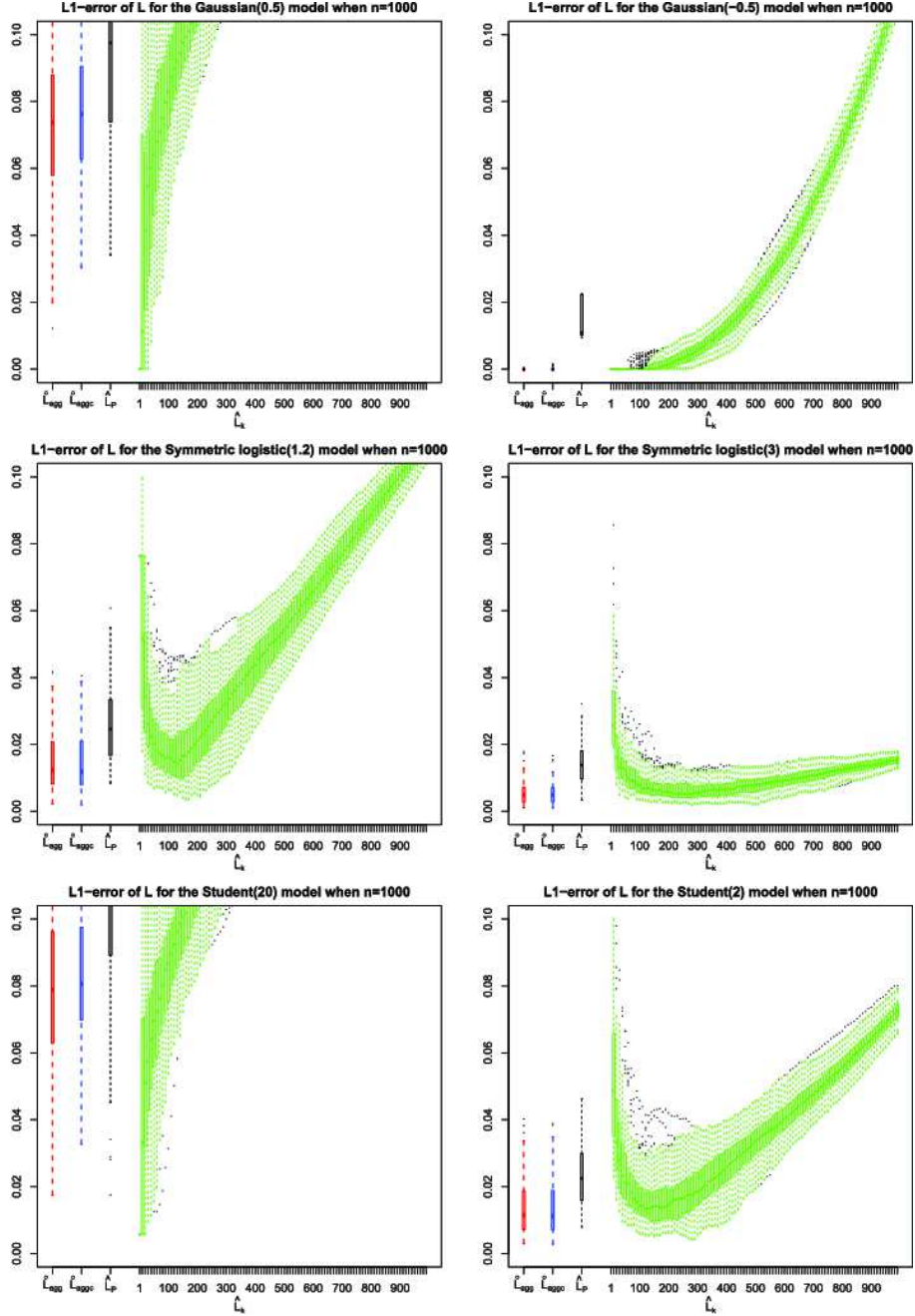}

\caption{Boxplot of the $L^1$-error of function $L$ for the estimators
$\protect\mathring L_{\mathrm{agg}}, \protect\mathring L_{\mathrm{aggc}}, \hat L_P$ and $\{
\hat L_k, k\}$.
First row: bivariate Normal model with correlation $\tau$: $\tau=0.5$
(left) and $\tau=-0.5$ (right).
Second row: bivariate Symmetric logistic$(s)$ model: $s=1/1.2$ (left)
and $s=1/3$ (right).
Third row: bivariate $\operatorname{Student}(\nu)$ model: $\nu=20$ (left) and $\nu
=2$ (right).}
\label{figl1-norm-L-rev1}
\end{figure}

Our results, summarized in Figure~\ref{figl1-norm-L-rev1}, will thus
exhibit in particular how the performance in the estimation of the s.t.d.f.
depends on the distance to the asymptotic independence case. The
$y$-axis scale has been fixed for all the six cases so that one can
measure that the estimation of the s.t.d.f. is a more ambitious problem
under asymptotic independence. However, our estimator $\mathring L_{\mathrm{agg}}$ has still nice properties when comparing it to the empirical
estimator $\hat{L}_k$.

The convex version $\mathring L_{\mathrm{aggc}}$ performs quite
equivalently as $\mathring L_{\mathrm{agg}}$. A reason for this is that by
construction our estimator is actually not far from a convex function.
So balancing the cost of convexifying with the benefit in the
performance motivates the simple use of $\mathring L_{\mathrm{agg}}$.

Finally, regarding Peng's estimator $\hat{L}_P$, one observes that
this estimator is an interesting alternative to the original family $\{
\hat{L}_k,k\}$, which, however, never outperforms our proposal.
%

\subsection{Estimating a failure probability}
{Let us illustrate in this subsection the question of estimating an
arbitrarily chosen failure probability $P(X>10^4$ or $Y>2
\cdot10^4)$, where $(X,Y)$ comes from the $\operatorname{BPII}(3)$ model, so that
$P(X>10^4$ or $Y>2 \cdot10^4) = 0.00011665$. Samples of size
$n=1000$ are considered. Thus empirical estimation} will be useless for
evaluating the probability of exceeding such extreme values for $X$ or
$Y$, and an extrapolation based on extreme value theory is thus
needed.

\begin{figure}

\includegraphics{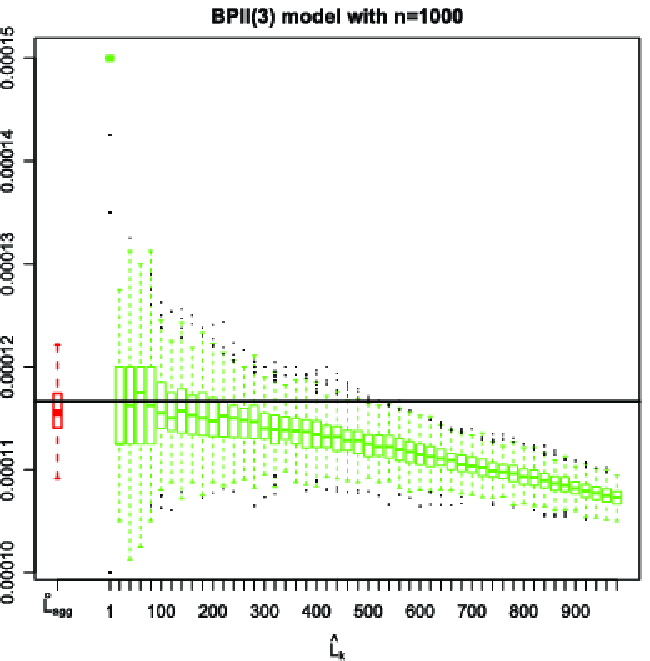}

\caption{Boxplot (based on 500 replicates) for the estimation of
$P(X>10^4$ or $Y>2\cdot10^4)$ when $(X,Y)$ is drawn from the
$\operatorname{BPII}(3)$ model with sample size $n=1000$ and assuming margins to be known.}\vspace*{-3pt}\label{figproba}
\end{figure}

First assume that it is known that the margins are standard Pareto.
This probability can be approximated by
\[
P\bigl(X>10^4 \mbox{ or } Y>2\cdot10^4\bigr)\simeq
\bigl(10^{-4}+5 \cdot 10^{-5} \bigr) L(2/3,1/3),
\]
which naturally comes from~(\ref{eqL}), the projection on the simplex
and the homogeneity of $L$. Estimating the unknown parameter
$L(2/3,1/3)$ with our candidate $\mathring L_{\mathrm{agg}}$ and the
original family $\{\hat L_k, k\}$ gives several boxplots (based on 500
replicates) that are presented in Figure~\ref{figproba}. The
comparison of these estimates is again favorable to~$\mathring L_{\mathrm{agg}}$.
%

\begin{figure}

\includegraphics{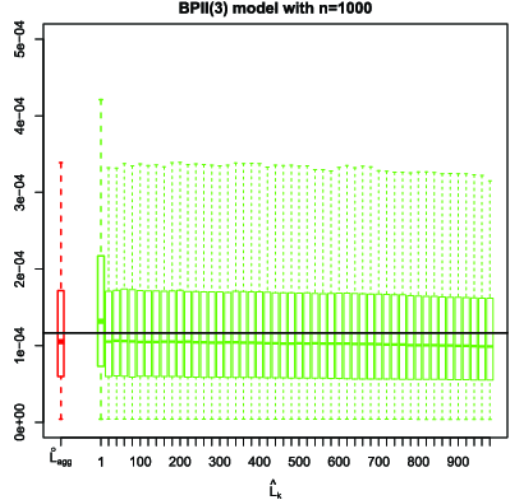}

\caption{Boxplot (500 replicates) of the estimation of $P(X>10^4$
or $Y>2\cdot10^4)$ when $(X,Y)$ is drawn from the $\operatorname{BPII}(3)$ model
with sample size $n=1000$ and estimating margins by POT method.}
\label{figprobabis}
\end{figure}

%
\begin{remark}
We also investigated the possible use of a second-order term in the
approximation of the probability
$P(X>10^4$ or $ Y>2\cdot10^4)$, making use of the following estimators
\[
\bigl(10^{-4}+5 \cdot10^{-5} \bigr) \mathring
L_{\mathrm{agg}} \biggl(\frac{2}{3},\frac{1}3 \biggr) + \biggl(
\frac{k}{n} \biggr)^{\hat
\rho} \bigl(10^{-4}+5
\cdot10^{-5} \bigr)^{1-\hat\rho} \hat \Delta_{k,2^{-1/\hat\rho}} \biggl(
\frac{2}{3},\frac{1}3 \biggr). %
\]
%
The results were so similar to those obtained in Figure~\ref
{figproba} that we chose to skip them.
\end{remark}
Second, when the margins are not assumed to be known, the estimation of
$p_1=1-F_1(10^4)$ and $p_2=1-F_2(2\cdot10^4)$ can be reached by the
POT method [see, e.g., \citet{beirlantgoegebeursegersteugels2004}, Section~7.4] for several
values of a threshold. After the study of mean residual life plots and
quantile plots, the thresholds have been fixed to be $X_{n-k,n}$ and
$Y_{n-k,n}$ for $k=200$.
The POT estimates deduced with these thresholds are, respectively,
denoted by $\hat{p}_1$ and $\hat{p}_2$.
The probability $P(X>10^4$ or $Y>2 \cdot10^4)$ is then
approximated by
\[
P\bigl(X>10^4\mbox{ or } Y>2\cdot10^4\bigr)\simeq (
\hat{p}_1+\hat {p}_2 ) L \biggl(\frac{\hat{p}_1}{\hat{p}_1+\hat{p}_2},
\frac
{\hat{p}_2}{\hat{p}_1+\hat{p}_2} \biggr).
\]
Estimating on each repetition the unknown parameter
$L (\hat{p}_1/(\hat{p}_1+\hat{p}_2),\hat{p}_2/(\hat{p}_1+\hat
{p}_2) )$
with our candidate $\mathring L_{\mathrm{agg}}$ and the original family $\{
\hat L_k, k\}$ gives several boxplots (based on 500 replicates)
presented in Figure~\ref{figprobabis}.
%
It seems clear that the uncertainty on the margins $F_1$ and $F_2$ has
much more influence than that of the s.t.d.f. $L$. Such findings
corroborate previous studies; see, for example, \citet{bruuntawn1998}
and \citet{dehaansinha1999}.\vadjust{\goodbreak}
\subsection{$Q$-curves}
\label{subsecQcurve}
\begin{figure}

\includegraphics{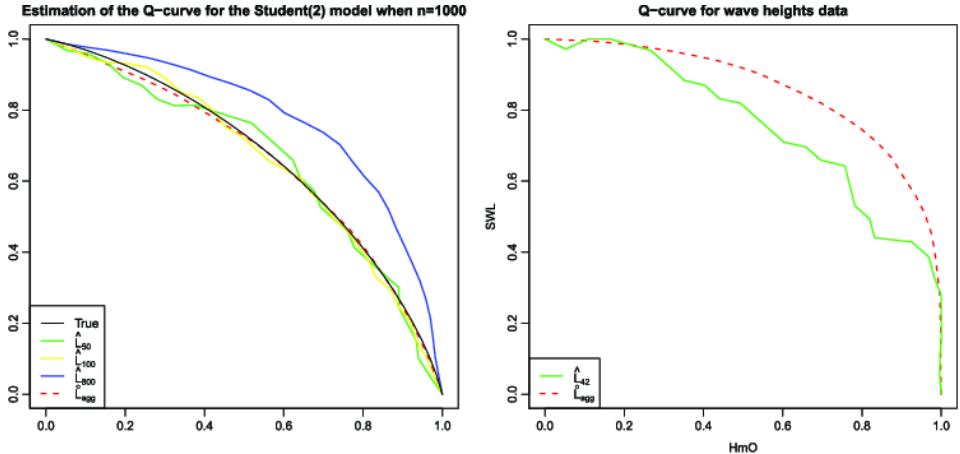}

\caption{Left: Estimation of the $Q$-curve for the bivariate
$\operatorname{Student}(2)$ law based on a sample of size 1000. Right: Estimated
$Q$-curve for {the wave heights data} introduced in \citet
{dehaanferreira2006}.}
\label{figneptune}
\end{figure}

\begin{figure}

\includegraphics{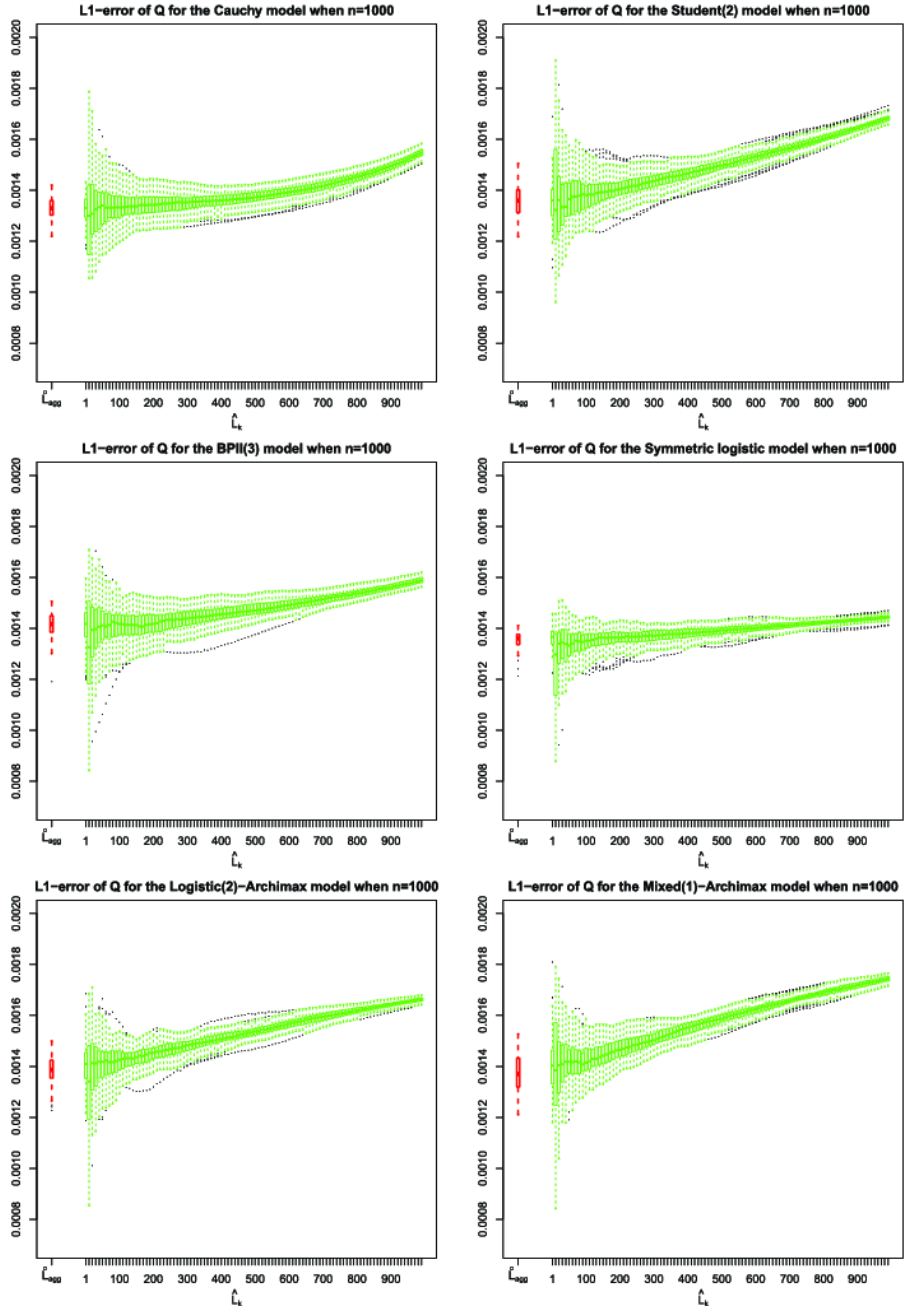}

\caption{Boxplot of the $L^1$-error of $Q$-curve for the estimators
$\protect\mathring L_{\mathrm{agg}}$ and $\{\hat L_k, k\}$.
First row: bivariate Cauchy model (left) and bivariate $\operatorname{Student}(2)$ model (right).
Second row: bivariate $\operatorname{BPII}(3)$ model (left) and bivariate Symmetric
logistic model (right).
Third row: bivariate Archimax model with logistic (left) and mixed
generator (right).}\label{figl1-norm-Q}
\end{figure}

Another nice representation of a function of several variables is
through its level sets. In the case of the function $L$, it consists of
looking (for any positive real $c$) at sets of the form $\{(x,y)\in
\mathbb{R}_+^2, L(x,y)\leq c\}$.
From homogeneity property, it is characterized by
\[
Q:=\bigl\{(x,y)\in\mathbb{R}_+^2, L(x,y)\leq1\bigr\}.
\]
Following \citeauthor{dehaanferreira2006} [(\citeyear{dehaanferreira2006}), page~245], the boundary of
this set can be written as
%
\[
\partial Q= \bigl\{ \bigl(b(\theta)\cos\theta,b(\theta)\sin\theta \bigr):
b(\theta)=\bigl(L(\cos\theta,\sin\theta)\bigr)^{-1}, \theta
\in[0,\pi/2] \bigr\}.
\]
The estimation of $\partial Q$ is naturally obtained by replacing $L$
by any estimator, and this is done here for the estimators $\mathring
L_{\mathrm{agg}}$ and $\{\hat L_k, k\}$.
Figure~\ref{figneptune} (left) exhibits the bias phenomenon (as $k$
increases) induced by $\hat L_k$ in the estimation of the $Q$-curve.
The bias factor on $\hat L_k$ is illustrated with $k=50, k=100$ and
$k=800$. The correction of the bias with $\mathring L_{\mathrm{agg}}$ is effective.
As in the previous section, the comparison of the different estimators
is provided in terms of a global criterium based on the $L^1$-norm,
given by
\[
\frac{\pi}{2(T+1)}\sum_{t=0}^T \biggl
\llvert \hat b \biggl(\frac{\pi
t}{2T} \biggr) - b \biggl(\frac{\pi t}{2T}
\biggr) \biggr\rrvert \biggl\{ \cos \biggl(\frac{\pi t}{2T} \biggr)+\sin \biggl(
\frac{\pi
t}{2T} \biggr) \biggr\}. %
\]
Figure~\ref{figl1-norm-Q} displays the boxplots of this measure,
based on $N=100$ realizations and for $T=30$ under the six bivariate
models given in the previous section.
%

The estimation of the $Q$-curve based on the original estimator $\hat
L_k$ is strongly sensitive to the choice of $k$: the bias (resp., the
variability) is an increasing (resp., decreasing) function of~$k$. The
performances of $\mathring L_{\mathrm{agg}}$ is similar to that of the
best~$\hat L_k$, which is unknown in practice. These features
corroborate the conclusions drawn in Section~\ref{subsecL1}.

To close this section, let us illustrate the $Q$-curve estimation {on
the wave heights data set} of \citet{dehaanferreira2006}, page~207. As explained therein, wave height (HmO) and
still water level (SWL) have been recorded during 828 storm events on
the Dutch coast. The analogous Figure~7.2 from \citet
{dehaanferreira2006} is reported in Figure~\ref{figneptune}
(right). Even if the two curves are not so close, the conclusion
remains the same: the estimated boundary is concave, which indicates
that the high values of the two variables are dependent.
%
%
\section{Estimation of second-order components \texorpdfstring{$\rho$}{rho} and $M$}\label{secrho}
In this section, we focus on the estimation of the function $M$ coming
from the second-order condition~(\ref{eq2ndorder}) and on the
estimation of its homogeneity parameter $1-\rho$.
\subsection{Second-order parameter \texorpdfstring{$\rho$}{rho}}
A possible way to estimate $\rho$ is to use
on each margin one of the techniques developed in the univariate
setting; see, for example, \citet{gomesdehaanpeng2002} or \citet
{ciupercamercadier2010}. Other methods make use of the multivariate
structure of the data; see, for example,
\citet{peng2010} and also \citet{goegebeurguillou2013} in a slightly
different framework. The construction described here takes likewise
advantage of the multivariate information of the sample.
With this purpose, the following proposition shows that a variable of
interest is the ratio of two terms $\hat\Delta_{k,a}$, defined
by~(\ref{eqDelta}).

\begin{proposition}
\label{proplim-quotient}
Assume that the conditions of Proposition~\ref{propcv-ps-L} are
fulfilled and fix positive real numbers $r$ and $a\in(0,1)$. Assume
moreover that the function $M$ never vanishes except on the axes.
Then, as $n$ tends to infinity, for every $\varepsilon>0$ and~$T>0$,
\[
\sup_{\varepsilon\leq x_1,\ldots,x_d \leq T} \biggl\llvert \frac{ \hat
\Delta_{k,a}(r {\mathbf x})} { \hat\Delta_{k,a}({\mathbf x}) } -
r^{1-\rho} \biggr\rrvert \stackrel{\mathbb{P}} {\longrightarrow} 0.
\]
\end{proposition}

%
\begin{remark} If the requirement that the function $M$ is either
positive, or negative in the positive quadrant does not hold, one could consider
the integral of $(\hat{\Delta}_{k,a}({\mathbf x}))^2$ over the set $\{
{\mathbf x}=(x_1,\ldots,x_d)$ s.t. $x_1^2+\cdots+x_d^2=1\}$ and prove a
result like Lemma~\ref{lemdelta} for this statistic. Then the
integral of $M^2$ appears in the denominator in Proposition~\ref
{proplim-quotient} 
instead of $M$ itself, and the sign of $M$ does not matter. This will
be part of a future work.
\end{remark}
%
A family of consistent estimators of the parameter $\rho$ can be
derived from Proposition~\ref{proplim-quotient}.
%
\begin{equation}
\label{defrho-hat} \hat\rho_{k,a,r}({\mathbf x}):= \biggl(1 -
\frac{1}{\log r} \log\biggl\llvert \frac{ \hat\Delta_{k,a}(r {\mathbf x})} { \hat\Delta_{k,a}({\mathbf x})
} \biggr\rrvert \biggr)
\wedge0.
\end{equation}
%
The following property can be obtained from the asymptotic expansion
given in Proposition~\ref{propdev-asympt-L}.

\begin{proposition}\label{prop3estim-rho}
Assume that the conditions of Proposition~\ref{propdev-asympt-L}
are fulfilled, and fix positive real numbers $r$ and $a\in(0,1)$.
Consider the estimator of $\rho$ defined by (\ref{defrho-hat}).
Assume moreover that the function $M$ never vanishes except on the axes.
Then, as $n$ tends to infinity,
\[
\sqrt{k} \alpha\biggl(\frac{n}{k}\biggr) \bigl\{ \hat
\rho_{k,a,r}({\mathbf x}) - \rho \bigr\} \overset{{d}} {\longrightarrow}\hat
Z_{ \rho, a, r}({\mathbf x}),
\]
in $D([\varepsilon,T]^d)$ for every $\varepsilon>0$ and $T>0$, with
\[
\hat Z_{ \rho, a, r}({\mathbf x}): = \frac{a^{-1} Z_L (a {\mathbf
x} )-Z_L({\mathbf x})}{ (a^{-\rho} -1) M({\mathbf x}) \log r } -\frac
{a^{-1} Z_L (ra {\mathbf x} )-Z_L(r {\mathbf x})}{ (a^{-\rho} -1)
M({\mathbf x}) r^{1-\rho} \log r}.
\]
\end{proposition}

\begin{figure}

\includegraphics{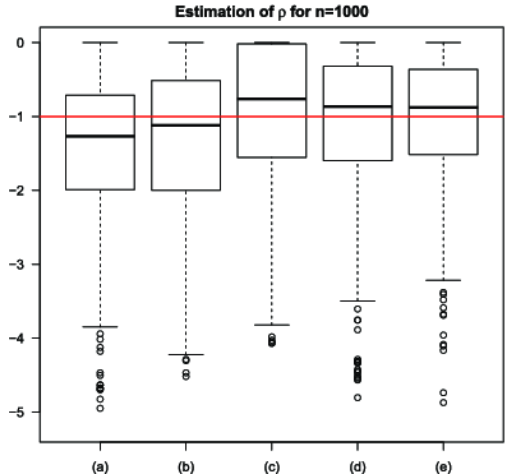}

\caption{Boxplot of 500 estimations of $\rho$ given by~(\protect\ref{defrho-hat}) using samples of size 1000 drawn from six models:
\textup{(a)}~$\operatorname{Student}(2)$;
\textup{(b)}~$\operatorname{BPII}(3)$;
\textup{(c)}~Symmetric Logistic with $s=1/3$;
\textup{(d)}~Archimax model with logistic generator with $s=1/2$;
\textup{(e)}~Archimax model with mixed generator.
Red line indicates the true value of $\rho=-1$.}
\label{figbox-rho}
\end{figure}

Figure~\ref{figbox-rho} illustrates the finite sample behavior of
this estimator of $\rho$ for a collection of bivariate models
introduced in Section~\ref{secexamples}, for which the true value of
$\rho$ is equal to $-$1.
%
These boxplots show that {the estimator performs reasonably well in
median, no matter the choice of model, but the uncertainty is rather
important. Fortunately this seems from simulation studies to have only
minor influence on the estimation of $L$.}
\subsection{Second-order function $M$}
Recall that from~(\ref{eqdev-asympt-La-L}) the asymptotic bias of
$\hat L_{k,a}({\mathbf x})$ is given by $\alpha(\frac{n}{k}) a^{-\rho}
M({\mathbf x})$. In order\vspace*{1pt} to circumvent an estimation of the term $\alpha
(n/k)$, a renormalization\vadjust{\goodbreak} is needed, focusing, for instance, on the
estimation of $M({\mathbf x})/M({\mathbf1/2})$ where ${\mathbf1/2}=(1/2,\ldots,1/2)$.
Thanks to~(\ref{eqDeltaCV}), this ratio can be consistently estimated by
\[
\frac{ \hat\Delta_{k,a}({\mathbf x})}{\hat\Delta_{k,a}({\mathbf1/2})}
\]
as soon as $k$ is a well-chosen intermediate sequence. The asymptotic
normality can also be derived from analogous arguments to those used in
the proof of Proposition~\ref{prop3estim-rho}. Details are not
presented here for the sake of simplicity.

\begin{figure}

\includegraphics{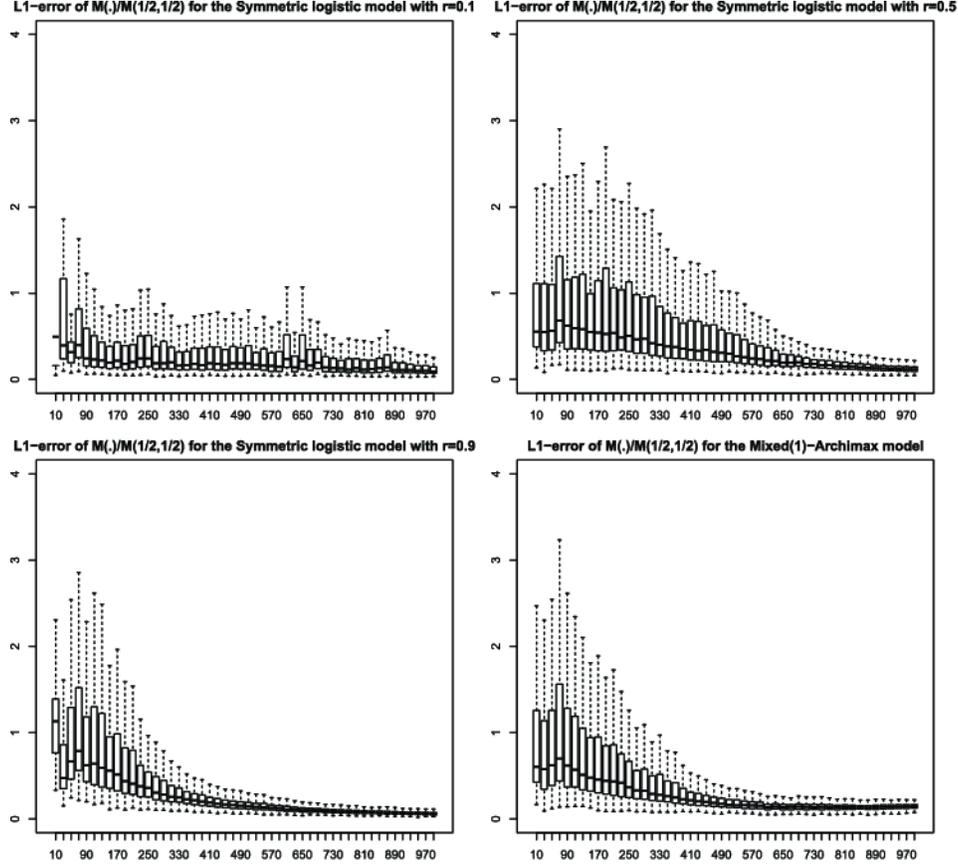}

\caption{Boxplot of the $L^1$-error of $M(\cdot)/M(1/2,1/2)$-curve.
First row: bivariate logistic model with $s=0.1$ (left) and with
$s=0.5$ (right).
Second row: bivariate logistic model with $s=0.9$ (left) and bivariate
Archimax with mixed generator (right).}\label{figl1-norm-M}
\end{figure}

Figure~\ref{figl1-norm-M} summarizes the behavior of the estimator of
the curve $t\mapsto M(t,1-t)/M(1/2,1/2)$ through boxplots of the
$L^1$-error, defined as in~(\ref{eqnorm1}). We observe from this
figure that the best estimation is reached for large values of $k$.
This feature does not depend on the degree of asymptotic dependence in
the Symmetric logistic model, nor on the strength of the bias of the
original estimator $\hat L_k$ detected on Figure~\ref{figl1-norm-L}.
These graphs confirm that the asymptotic bias is remarkably well
estimated for large values of $k$. This helps to understand why the
bias subtraction is accurate for large or very large choices of $k$, as
also commented in Section~\ref{subseccorrections}.
%
\section{Concluding comments}

This paper deals with the estimation of the extremal dependence
structure in a multivariate context. Focusing on the s.t.d.f., the
empirical counterpart is the nonparametric reference. A common feature
when modeling extreme events is the delicate choice of the number of
observations used in the estimation, and it spoils the good performance
of this estimator. The aim of this paper has been to correct the
asymptotic bias of the empirical estimator, so that the choice of the
threshold becomes less sensitive. Two asymptotically unbiased
estimators have been proposed and studied, both theoretically and
numerically. The estimator defined in Section~\ref{subsecmethodB}
proves to outperform the original estimator, whatever the model
considered. Its aggregated version defined in Section~\ref
{subseccorrections} appears as a worthy candidate to estimate the s.t.d.f.


\section{Proofs}\label{secproofs}

\mbox{}
\begin{pf*}{Proof of Proposition~\ref{propcv-ps-L}}
Denote by $U^{(j)}_i$ the uniform random variables
$U^{(j)}_i=1-F_j(X^{(j)}_i)$ for $j=1,\ldots,d$.
Introducing
\[
V_k({\mathbf x})= \frac{1}k \sum_{i=1}^n
\ind_{ \{U^{(1)}_i \leq
kx_1/n~\mathrm{or}~\ldots~\mathrm{or}~U^{(d)}_i \leq kx_d/n
 \} }
\]
allows us to rewrite $\hat L_k$ as the following:
\[
\hat L_k ({\mathbf x}) = V_k \biggl(\frac{n}{k}
U^{(1)}_{[kx_1],n}, \ldots, \frac{n}{k}U^{(d)}_{[kx_d],n}
\biggr).
\]
Write
\begin{eqnarray*}
&& \hat L_k({\mathbf x}) - L({\mathbf x})
\\
&&\qquad = V_k
\biggl(\frac{n}{k} U^{(1)}_{[kx_1],n}, \ldots,
\frac{n}{k}U^{(d)}_{[kx_d],n} \biggr)
\\
&&\quad\qquad{}  - \frac{n}{k}
\bigl[ 1- F\bigl\{F_1^{-1}\bigl(1-U^{(1)}_{[kx_1],n}
\bigr),\ldots,F_d^{-1}\bigl(1-U^{(d)}_{[kx_d],n}
\bigr) \bigr\}\bigr]
\\
&&\quad\qquad{}+ \frac{n}{k} \bigl[ 1- F\bigl\{F_1^{-1}
\bigl(1-U^{(1)}_{[kx_1],n}\bigr),\ldots,F_d^{-1}
\bigl(1-U^{(d)}_{[kx_d],n}\bigr) \bigr\}\bigr]
\\
&&\quad\qquad{} - L \biggl(
\frac{n}{k} U^{(1)}_{[kx_1],n}, \ldots, \frac{n}{k}U^{(d)}_{[kx_d],n}
\biggr)
\\
&&\quad\qquad{} + L \biggl(\frac{n}{k} U^{(1)}_{[kx_1],n}, \ldots,
\frac
{n}{k}U^{(d)}_{[kx_d],n} \biggr) - L({\mathbf x}),
\end{eqnarray*}
and denote $A_{1,k}({\mathbf x})$ [resp., $A_{2,k}({\mathbf x})$ and
$A_{3,k}({\mathbf x})$] the first line (resp., second and third lines) of
the right-hand side.

Applying \citeauthor{dehaanferreira2006} [(\citeyear{dehaanferreira2006}), Proposition 7.2.3] leads to
\[
\sqrt{k} A_{1,k}({\mathbf x}) \stackrel{d} {\to} W_L({\mathbf x}),
\]
in $D([0,T]^d)$ for every $T>0$ and for any intermediate sequence,
where $W_L$ is a continuous centered Gaussian process with covariance
structure specified in
Proposition~\ref{propdev-asympt-L}.
Due to the Skorohod construction we can write
%
\begin{equation}
\label{eqskorohodA1} \sup_{0\leq x_1,\ldots, x_d \leq T} \bigl\llvert \sqrt{k}
A_{1,k}({\mathbf x}) - W_L({\mathbf x}) \bigr\rrvert \to0\qquad
\mbox{a.s.},
\end{equation}
which implies, since $\sqrt{k}\alpha(n/k)\to\infty$,
\[
\sup_{0\leq x_1,\ldots, x_d \leq T} \biggl\llvert \frac{A_{1,k}({\mathbf
x})}{\alpha (n/k)} \biggr\rrvert
=O_{\mathbb{P}} \biggl(\frac{1}{\sqrt{k}\alpha (n/k)} \biggr).
\]
Again for any intermediate sequence, the proof of \citeauthor{dehaanferreira2006}
[(\citeyear{dehaanferreira2006}), Theorem~7.2.2] ensures the convergence for $j=1,\ldots,d$
%
\begin{equation}
\label{eqEKSmarg1} \sup_{x\in[0,T]} \biggl\llvert \sqrt{k} \biggl(
\frac{n}{k}U^{(j)}_{[kx],n} -x \biggr) +
W_L(x{\mathbf e}_j)\biggr\rrvert \to0 \qquad\mbox{a.s.},
\end{equation}
and finally
%
\begin{equation}
\label{eqskorohodA3} \sup_{0\leq x_1,\ldots, x_d \leq T} \Biggl\llvert \sqrt{k}
A_{3,k}({\mathbf x}) +\sum_{j=1}^dW_L(x_j{
\mathbf e}_j) \partial_j L({\mathbf x}) \Biggr\rrvert \to 0
\qquad\mbox{a.s.}
\end{equation}
As previously, this yields
\[
\sup_{0\leq x_1,\ldots, x_d \leq T} \biggl\llvert \frac{A_{3,k}({\mathbf
x})}{\alpha (n/k)} \biggr\rrvert =O
\biggl(\frac
{1}{\sqrt{k}\alpha (n/k)} \biggr).
\]
Since the intermediate sequence satisfies $\sqrt{k}\alpha (\frac
{n}{k} ) \to\infty$, it thus remains to prove that
\[
\sup_{0\leq x_1,\ldots, x_d \leq T} \biggl\llvert \frac{A_{2,k}({\mathbf
x})}{\alpha(n/k)} -M({\mathbf x}) \biggr
\rrvert \to0 \qquad\mbox {a.s.}
\]
The second-order condition that holds uniformly on $[0,T]^d$ in~(\ref
{eq2ndorder}) yields
\[
\sup_{0\leq x_1,\ldots, x_d \leq T} \biggl\llvert \frac{A_{2,k}({\mathbf
x})}{\alpha(n/k)} -M \biggl(
\frac{n}{k} U^{(1)}_{[kx_1],n}, \ldots, \frac{n}{k}U^{(d)}_{[kx_d],n}
\biggr)\biggr\rrvert \to0 \qquad\mbox{a.s.} %
\]
Then the result follows from
\[
\sup_{0\leq x_1,\ldots, x_d \leq T} \biggl\llvert M({\mathbf x}) -M \biggl(
\frac
{n}{k} U^{(1)}_{[kx_1],n}, \ldots, \frac{n}{k}U^{(d)}_{[kx_d],n}
\biggr)\biggr\rrvert \to0 \qquad\mbox{a.s.}, %
\]
which is obtained combining~(\ref{eqEKSmarg1}) and the continuity of
the function $M$.
\end{pf*}

\begin{pf*}{Proof of Proposition~\ref{propdev-asympt-L}}
We use the notation introduced in the proof of Proposition~\ref
{propcv-ps-L}. Thanks to the Skorohod construction, we can start
from~(\ref{eqskorohodA1}). Combined with~(\ref{eqskorohodA3}), it
is sufficient to prove the convergence
\[
\sup_{0\leq x_1,\ldots, x_d \leq T} \biggl\llvert \sqrt{k} \biggl\{ A_{2,k}({
\mathbf x}) -\alpha\biggl(\frac{n}{k}\biggr) M({\mathbf x}) \biggr\} \biggr\rrvert
\to0 \qquad\mbox{a.s.}
\]
Note that the third-order condition, the uniformity on $[0,T]^d$ of the
convergence in~(\ref{eq3rdorder}) and the continuity of $N$ yield
\[
A_{2,k}({\mathbf x}) = \alpha\biggl(\frac{n}{k}\biggr)M \biggl(
\frac{n}{k} U^{(1)}_{[kx_1],n}, \ldots, \frac{n}{k}U^{(d)}_{[kx_d],n}
\biggr)+ O_{\mathbb P} \biggl(\alpha\biggl(\frac{n}{k}\biggr)\beta
\biggl(\frac{n}{k}\biggr) \biggr). %
\]
Thanks to~(\ref{eqEKSmarg1}) and to the existence of the first-order
partial derivatives $\partial_j M$ $(j = 1, \dots, d)$ of the
function $M$, we have that
\begin{eqnarray*}\label{eqdev-M}
&& \sup_{0\leq x_1,\ldots, x_d \leq T} \Biggl\llvert \sqrt{k} \biggl\{ M
\biggl(\frac{n}{k} U^{(1)}_{[kx_1],n}, \ldots,
\frac{n}{k}U^{(d)}_{[kx_d],n} \biggr) - M({\mathbf x}) \biggr\}
\\
&&\hspace*{121pt}\qquad{}  +
\sum_{j=1}^dW_L(x_j{
\mathbf e}_j) \partial_j M({\mathbf x}) \Biggr\rrvert
\end{eqnarray*}
converges to 0 in probability, as $n$ tends to infinity. This implies that
\[
\sup_{0\leq x_1,\ldots, x_d \leq T} \biggl\llvert \sqrt{k} \biggl\{ A_{2,k}({
\mathbf x}) -\alpha\biggl(\frac{n}{k}\biggr) M({\mathbf x}) \biggr\} \biggr\rrvert
= O_{\mathbb P} \biggl(\biggl\llvert \sqrt{k}\alpha\biggl(\frac{n}{k}
\biggr)\beta\biggl(\frac{n}{k}\biggr) + \alpha\biggl(\frac{n}{k}
\biggr)\biggr\rrvert \biggr),
\]
which completes the proof, thanks to the choice of the intermediate sequence.\vadjust{\goodbreak}
\end{pf*}

\begin{pf*}{Proof of Theorem~\ref{thmbiasB}}
Recall that $b=(a^{-\rho} +1)^{-1/\rho}$, and denote $\hat
b=(a^{-\hat\rho} +1)^{-1/\hat\rho}$. Write
%
\begin{equation}
\mathring L_{k, a, k_\rho}- L= \{ \hat L_{k,a} - L \} + \{ \hat
L_k - L\} - \{ \hat L_{k, \hat b} - L\},\label{eqdecomp1}
\end{equation}
%
which equals, thanks to~(\ref{eqdev-asympt-La-L}) and under
Skorohod's construction,
\begin{eqnarray*}
&& \alpha \biggl(\frac{n}{k} \biggr) \bigl(a^{-\rho} + 1\bigr) M({
\mathbf x}) + \frac{1}{\sqrt k} \bigl(a^{-1} Z_L(a {\mathbf x}) +
Z_L({\mathbf x}) \bigr)
\\[-1pt]
&&\quad{} - \alpha \biggl(\frac{n}{k} \biggr) \hat
b^{-\rho} M({\mathbf x}) - \frac
{ b^{-1}}{\sqrt k} Z_L(b {\mathbf x})+ o
\biggl(\frac{1}{\sqrt{k}} \biggr)
\\[-1pt]
&&\qquad  = \alpha \biggl(\frac{n}{k} \biggr) \bigl( \bigl(a^{-\rho} + 1
\bigr)-b^{-\rho} \bigr) M({\mathbf x})+ \frac{1}{\sqrt k}\mathring
Y_a({\mathbf x})
\\[-1pt]
&&\quad\qquad{} +\alpha \biggl(\frac{n}{k} \biggr)
\bigl(b^{-\rho}-\hat b^{-\rho
} \bigr)M({\mathbf x})+ o \biggl(
\frac{1}{\sqrt{k}} \biggr)
\\[-1pt]
&&\qquad = \alpha \biggl(\frac{n}{k} \biggr) \bigl( \bigl(a^{-\rho} + 1
\bigr)-b^{-\rho} \bigr) M({\mathbf x})+ \frac{1}{\sqrt k}\mathring
Y_a({\mathbf x})
\\[-1pt]
&&\quad\qquad{} +\alpha \biggl(\frac{n}{k} \biggr)O_\mathbb{P}
\biggl(\frac
{1}{\sqrt{k_\rho} \alpha(n/k_\rho)} \biggr)+ o \biggl(\frac
{1}{\sqrt{k}} \biggr).
\end{eqnarray*}
The first term is zero. Since both $k = o(k_\rho)$ and $\alpha$ is
regularly varying with negative index, the only the last term can be
put into the term $o (\frac{1}{\sqrt{k}} )$.
Finally, the covariance function follows from the equality in law as
processes between $Z_L(a {\mathbf x}) $ and $\sqrt{a} Z_L({\mathbf x})$.
\end{pf*}

The proofs of Theorem~\ref{thmbias-Lu} and Proposition~\ref
{prop3estim-rho} are based on the following auxiliary result.

\begin{lemma}\label{lemdelta}
Assume that the conditions of Proposition~\ref{propdev-asympt-L} are
fulfilled. Then for any positive real $r$, one has as $n$ tends to infinity,
\begin{eqnarray*}
&& \sqrt{k} \alpha\biggl(\frac{n}{k}\biggr) \biggl\{ \frac{\hat\Delta_{k,a}(r
{\mathbf x})}{\alpha(n/k)} -
\bigl(a^{-\rho} -1\bigr)r^{1-\rho} M({\mathbf x}) \biggr\}
\stackrel{d} {
\to} a^{-1} Z_L (ra {\mathbf x} ) - Z_L(r {\mathbf x}),
\end{eqnarray*}
in $D([0,T]^d)$ for every $T>0$.
\end{lemma}

\begin{pf*}{Proof of Lemma \ref{lemdelta}}
Making use of the homogeneity of the function $L$, write
\[
\hat\Delta_{k,a}(r {\mathbf x}) = \bigl\{ \hat L_{k,a}(r {\mathbf x})
- L(r {\mathbf x}) \bigr\} - \bigl\{ \hat L_k(r {\mathbf x}) -L(r {\mathbf x})
\bigr\}. %
\]
Using the Skorohod construction, it follows from equations~(\ref
{eqasympt-dev-L}) and~(\ref{eqdev-asympt-La-L}) that
\begin{eqnarray*}
&& \sup_{0\leq x_1,\ldots, x_d \leq T/r} \biggl\llvert \sqrt{k} \alpha\biggl(
\frac{n}{k}\biggr) \biggl\{\frac{\hat\Delta_{k,a}(r {\mathbf x})}{\alpha(
n/k)}-\bigl(a^{-\rho} -1 \bigr)r^{1-\rho} M({\mathbf x}) \biggr\}
\\
&&\hspace*{149pt}{} -a^{-1} Z_L (ra {
\mathbf x} ) + Z_L(r {\mathbf x}) \biggr\rrvert
\end{eqnarray*}
tends to 0 almost surely, as $n$ tends to infinity.
\end{pf*}

\begin{pf*}{Proof of Theorem~\ref{thmbias-Lu}}
Note that
\begin{eqnarray*}
&& \hat L_k({\mathbf x}) \frac{\hat\Delta_{k_\rho,a}(a {\mathbf x})}{\alpha
(n/k_\rho)} -\hat L_k(a {
\mathbf x}) \frac{ \hat\Delta_{k_\rho,a}({\mathbf
x})}{\alpha(n/k_\rho)}
\\
&&\qquad  = \hat L_k({\mathbf x}) \biggl(
\frac{\hat\Delta
_{k_\rho,a}(a {\mathbf x})}{\alpha(n/k_\rho)} - a \frac{ \hat\Delta
_{k_\rho,a}({\mathbf x})}{\alpha(n/k_\rho)} \biggr) -a\frac{\hat\Delta
_{k_\rho,a}({\mathbf x})\hat\Delta_{k,a}({\mathbf x})}{\alpha(n/k_\rho)}.
\end{eqnarray*}
Under a Skorohod construction, Lemma~\ref{lemdelta} allows us to
write the expansions of the terms $\hat\Delta_{k,a}({\mathbf x})$, $\hat
\Delta_{k_\rho,a}({\mathbf x})$ and $\hat\Delta_{k_\rho,a}(a {\mathbf
x})$, which implies on the one hand
%
\begin{eqnarray}\label{eqdeltaterm1}
&& \frac{\hat\Delta_{k_\rho,a}(a {\mathbf x})}{\alpha(n/k_\rho)} - a \frac{ \hat\Delta_{k_\rho,a}({\mathbf x})}{\alpha(n/k_\rho)}\nonumber
\\
&&\qquad  =a \bigl(a^{-\rho} -1 \bigr)^2 M({\mathbf x})
\nonumber\\[-8pt]\\[-8pt]\nonumber
&&\qquad\quad {} + \frac{1}{\sqrt{k_\rho} \alpha(n/k_\rho)} \bigl\{ a^{-1} Z_L
\bigl(a^2 {\mathbf x}\bigr) -2 Z_L(a {\mathbf x}) + a
Z_L({\mathbf x}) \bigr\}
\\
&&\quad\qquad{}  +o \biggl( \frac{1}{\sqrt{k_\rho} \alpha(n/k_\rho)} \biggr),\nonumber
\end{eqnarray}
and
%
\begin{eqnarray}\label{eqdeltaterm2}
\frac{\hat\Delta_{k_\rho,a}({\mathbf x})\hat\Delta_{k,a}({\mathbf
x})}{\alpha(n/k_\rho)} &=& \alpha(n/k) \bigl(a^{-\rho}-1\bigr)^2
M^2({\mathbf {x}})\nonumber
\\
&&{} +\bigl(a^{-\rho}-1\bigr)M({\mathbf x})
\frac{a^{-1}Z_L(a{\mathbf x})-Z_L({\mathbf
x})}{\sqrt{k}}
\\
&&{} + O_{\mathbb P} \biggl(\frac{\alpha(n/k)}{\sqrt{k_\rho} \alpha
(n/k_\rho)} + \frac{1}{\sqrt{k} \sqrt{k_\rho} \alpha(n/k_\rho
)} \biggr)+o \biggl(\frac{1}{\sqrt{k}} \biggr)\nonumber
\end{eqnarray}
on the other hand, both uniformly for ${\mathbf x}\in[\varepsilon,T]^d$.
Combining~(\ref{eqdeltaterm1}) and~(\ref{eqdeltaterm2}) with
equation~(\ref{eqasympt-dev-L}), one gets
\begin{eqnarray*}
&& \hat L_k({\mathbf x}) \frac{\hat\Delta_{k_\rho,a}(a {\mathbf x})}{\alpha
(n/k_\rho)} -\hat L_k(a {
\mathbf x}) \frac{ \hat\Delta_{k_\rho,a}({\mathbf
x})}{\alpha(n/k_\rho)}\nonumber
\\
&& \qquad =a \bigl(a^{-\rho} -1\bigr)^2 M({\mathbf x}) L({\mathbf x}) +
\frac{1}{\sqrt{k}} M({\mathbf x}) \bigl(a^{-\rho} -1\bigr)
\bigl(a^{1-\rho} Z_L({\mathbf x}) - Z_L(a {\mathbf x})
\bigr)
\\
&&\quad\qquad {}+ \frac{1}{\sqrt{k_\rho} \alpha(n/k_\rho)} L({\mathbf x}) \bigl\{ a^{-1} Z_L
\bigl(a^2 {\mathbf x}\bigr) -2 Z_L(a {\mathbf x}) + a
Z_L({\mathbf x}) \bigr\}
\\
&&\quad\qquad{} +o \biggl(\frac{1}{\sqrt{k}} \biggr) +o \biggl(
\frac{1}{\sqrt{k_\rho} \alpha(n/k_\rho)} \biggr).\nonumber
\end{eqnarray*}

Since the last expression and equation~(\ref{eqdeltaterm1}) are,
respectively, the numerator and denominator of $\tilde L_{k,k_\rho,
a}({\mathbf x})$, one obtains, after simplification,
\[
\sqrt{k} \bigl(\tilde L_{k,k_\rho, a}({\mathbf x}) - L({\mathbf x})\bigr) =
\frac
{a^{-\rho} Z_L({\mathbf x}) - a^{-1} Z_L(a {\mathbf x}) }{ a^{-\rho} -1} + o \biggl(\frac{\sqrt{k}}{ \sqrt{k_\rho} \alpha(n/k_\rho)} \biggr) +o(1),
\]
since $M$ does not vanish by assumption.
The choice of the sequences $k$ and $k_\rho$ allows us to conclude
since $\sqrt{k} = O (\sqrt{k_\rho} \alpha(n/k_\rho) )$.
\end{pf*}

\begin{pf*}{Proof of Proposition~\ref{proplim-quotient}}
Applying Lemma~\ref{lemdelta}, we have
%
\begin{equation}
\label{eqsupdelta} \sup_{\varepsilon\leq x_1,\ldots,x_d \leq T} \biggl\llvert \frac{\hat\Delta
_{k,a}({\mathbf x})}{\alpha (n/k)}-
\bigl(a^{-\rho
}-1\bigr)M({\mathbf x})\biggr\rrvert \stackrel{\mathbb{P}} {
\longrightarrow} 0.
\end{equation}
As a consequence,
\begin{eqnarray*}
&& \sup_{\varepsilon\leq x_1,\ldots,x_d \leq T} \biggl\llvert \frac{ \hat
\Delta_{k,a}(r {\mathbf x})} { \hat\Delta_{k,a}({\mathbf x}) } -
r^{1-\rho} \biggr\rrvert
\\
&&\qquad = \sup_{\varepsilon\leq x_1,\ldots,x_d \leq T} \biggl\llvert
\frac
{ \hat\Delta_{k,a}(r {\mathbf x})/\alpha(n/k)} { \hat\Delta_{k,a}({\mathbf
x})/\alpha(n/k) } - r^{1-\rho} \biggr\rrvert
\\
&&\qquad = O_{\mathbb P} \biggl( \sup_{\varepsilon\leq x_1,\ldots,x_d \leq T} \biggl\llvert
\frac{\hat\Delta_{k,a}(r {\mathbf x})}{\alpha(n/k)} - r^{1-\rho
} \frac{ \hat\Delta_{k,a}({\mathbf x})}{\alpha(n/k)} \biggr\rrvert \biggr),
\end{eqnarray*}
since $(a^{-\rho} -1)M({\mathbf x}) \neq0$ by assumption.
Writing
\begin{eqnarray*}
&& \biggl\llvert \frac{\hat\Delta_{k,a}(r {\mathbf x})}{\alpha(n/k)} - r^{1-\rho
} \frac{ \hat\Delta_{k,a}({\mathbf x})}{\alpha(n/k)} \biggr
\rrvert
\\
&&\qquad  \leq \biggl\llvert \frac{\hat\Delta_{k,a}(r {\mathbf x})}{\alpha(n/k)} - r^{1-\rho
}
\bigl(a^{-\rho}-1\bigr)M({\mathbf x}) \biggr\rrvert
\\
&&\quad\qquad{} + \biggl\llvert r^{1-\rho} \bigl(a^{-\rho}-1\bigr)M({\mathbf x}) -
r^{1-\rho} \frac{ \hat
\Delta_{k,a}({\mathbf x})}{\alpha(n/k)} \biggr\rrvert,
\end{eqnarray*}
and using twice equation~(\ref{eqsupdelta}) leads to the conclusion.
\end{pf*}

\begin{pf*}{Proof of Proposition~\ref{prop3estim-rho}}
Define $  Q_{k,a,r}({\mathbf x}):=\frac{ \hat\Delta_{k,a}(r
{\mathbf x})} { \hat\Delta_{k,a}({\mathbf x}) }$.
Lemma~\ref{lemdelta} used twice yields
%
\begin{equation}
\label{eqdev-asympt-Q} \sqrt{k} \alpha\biggl(\frac{n}{k}\biggr)
\bigl(Q_{k,a,r}({\mathbf x}) - r^{1-\rho}\bigr) \stackrel{d} {\to} -
r^{1-\rho} \log r \hat Z_{\rho,a, r}({\mathbf x}),
\end{equation}
where $ \hat Z_{\rho,a, r}({\mathbf x})$ is defined in Proposition~\ref
{prop3estim-rho}.
Since $\hat\rho_{k,a,r}({\mathbf x})=1-\log(Q_{k,a,r}({\mathbf x}))/ \log r$,
the result follows straightforwardly from~(\ref{eqdev-asympt-Q}) and
the Delta method.
\end{pf*}

\section*{Acknowledgments}

We wish to thank Armelle Guillou for pointing out a deficiency in the
original version of the paper, as well as several misprints. We thank
the referees for very helpful comments.




%

\printaddresses
\end{document}